%%%%%%%%%%below I am commenting out \input usual and inputting
%%%%%%%%%%usual.tex by hand. this is for submission to the archives.

 \overfullrule=0mm
\baselineskip=0,2cm
\lineskiplimit=1pt
\lineskip=1mm
\def \Ci {C^\infty}

\def\noi{{\noindent}}
\def \l {{\lambda}}
\def\R{{\bf R}}
\def\N{{\bf N}}
\def\entiers{{\rm I\kern-2pt N}}
\def\Z{{\bf Z}}
\def\C{{\bf C}}
\def\P{{\cal P}}
\def\G{{\cal G}}

\hsize 16 truecm
\vsize 24 truecm

\def\sp{{\vskip 0,2cm}}
\def\bp{{\vskip 0,5cm}}
\def\np{{\vskip 0,3cm}}
 \def \e {\varepsilon}

\def\Det{{\rm Det}}
\def\str{{\rm str}}
\def\tr{{\rm tr}}
\def\Ker{{\rm Ker}}
\def\ord{{\rm ord}}
\def\Coker{{\rm Coker}}
\def\Hom{{\rm Hom}}
\def\Lim{{\rm Lim}}
\def\res{{\rm res}}
\def\ch{{\rm ch}}
\def\op{{\rm op}}
\def\LC{{\rm LC}}
\def\Lie{{\rm Lie}}

\centerline{\bf CURVATURE ON DETERMINANT BUNDLES AND FIRST CHERN FORMS}
 \vskip 0,5cm
  \centerline{     Sylvie PAYCHA*, Steven ROSENBERG**  }
\vskip 0,3cm
 \centerline{(*)Laboratoire de Math\'ematiques Appliqu\'ees}
\centerline{Universit\'e Blaise Pascal (Clermont II)}
\centerline{Complexe Universitaire des C\'ezeaux}
\centerline{63177 Aubi\`ere Cedex}
  \vskip 0,3cm
\centerline{(**) Mathematics Department}
\centerline{Boston University}
\centerline{Boston, Massachusetts}
\vskip 0,5cm
 \centerline{\tt paycha@ucfma.univ-bpclermont.fr}
\centerline{\tt sr@math.bu.edu }
\vskip 2cm
 \noi
\centerline{ \bf Acknowledgements}
\np
The first author  would like to thank the Department of Mathematics 
and Statistics at Boston University
for its hospitality,
and 
the second author   expresses his thanks to the Mathematics Department of the
Universit\'e Blaise Pascal for its hospitality.
 \vskip 2cm
 \noi
\centerline{\bf  Abstract}\np

The Quillen-Bismut-Freed construction associates a determinant line
bundle with connection to an infinite dimensional super vector bundle with a
family of Dirac-type operators. We define the regularized
first Chern form of the infinite dimensional bundle, and relate it
to
the curvature of the Bismut-Freed
connection on the determinant bundle. In finite dimensions, these
forms agree (up to sign), but in infinite dimensions
there is a correction term, which
 we express in terms of Wodzicki residues.

\sp We illustrate these results with a
string theory computation.  There is
a natural super
vector bundle over the manifold
of smooth
almost complex structures
on a Riemannian surface.  The Bismut-Freed superconnection is identified
with classical Teichm\"uller theory connections, and its curvature and
regularized first Chern form are computed.
    \vfill 


\eject \noindent

%%%%%%%%%%%%%%%%%%%%%%%%%usual.tex%%%%%%%%%%%%%%%%%%%%%%%%
% This file defines the standard font setup, including the macros
% \tenpoint and \twelvepoint, as well as some standard dimension settings.
%
% It is indended to be used as a basis for most other ``personal'' formats.

% Make sure we haven't already been run

\ifx\UsualIsLoaded\undefined
\let\UsualIsLoaded=\relax		% define UsualIsLoaded
\newskip\ttglue

\font\twelverm=cmr12
\font\twelvebf=cmbx12
\font\twelvett=cmtt12
\font\twelvesl=cmsl12
\font\twelvesy=cmsy10 scaled \magstep1
\font\twelvei=cmmi12
\font\twelveit=cmti12
\font\twelvesc=cmcsc10 scaled \magstep1
\font\twelveBbb=bbmbx12 
\font\tenrm=cmr10
\font\tenbf=cmb10 
\font\tentt=cmtt10 
\font\tensl=cmsl10 
\font\tensy=cmsy10 
\font\teni=cmmi10 
\font\tenit=cmti10 
\font\tensc=cmcsc10 
\font\tenBbb=bbmbx10

%%%%
\font\elevenrm=cmr10 scaled 1100
\font\elevenbf=cmb10 scaled 1100
\font\eleventt=cmtt10 scaled 1100
\font\elevensl=cmsl10 scaled 1100
\font\elevensy=cmsy10 scaled 1100
\font\eleveni=cmmi10 scaled 1100
\font\elevenit=cmti10 scaled 1100

\font\elevenBbb=bbmbx10 scaled 1100
%%%%%
%
\font\ninei=cmmi9 
\font\ninerm=cmr9

\font\ninesy=cmsy9 
\font\nineBbb=bbmbx9
\font\eighti=cmmi8
\font\eightrm=cmr8

\font\eightsy=cmsy8
\font\eightBbb=bbmbx8
\font\seveni=cmmi7 
\font\sevenrm=cmr7
\font\sevensy=cmsy7
\font\sevenBbb=bbmbx7
%
% large fonts
%

%

%%%%%%%%%% define the BlackBoard family 
\newfam\Bbbfam

%%%%%%%%%%%%%% this macro defines a tenpoint font
%%%%%%
\def\tenpoint{
    \def\rm{\fam0\tenrm}
    \textfont0=\tenrm \scriptfont0=\eightrm \scriptscriptfont0=\sevenrm
    \textfont1=\teni \scriptfont1=\eighti \scriptscriptfont1=\seveni
    \textfont2=\tensy \scriptfont2=\eightsy \scriptscriptfont2=\sevensy
    \textfont3=\tenex \scriptfont3=\tenex \scriptscriptfont3=\tenex

    \textfont\itfam=\tenit \def\it{\fam\itfam\tenit}
    \textfont\slfam=\tensl \def\sl{\fam\slfam\tensl}
    \textfont\bffam=\tenbf \def\bf{\fam\bffam\tenbf}
    \textfont\ttfam=\tentt \def\tt{\fam\ttfam\tentt}
    \textfont\Bbbfam=\tenBbb \def\Bbb{\fam\Bbbfam\tenBbb}
    \scriptfont\Bbbfam=\eightBbb \scriptscriptfont\Bbbfam=\sevenBbb
    \def\sc{\tensc}
    \tt \ttglue=.5em plus.25em minus.15em
    \normalbaselineskip=12pt

    \setbox\strutbox=\hbox{\vrule height10pt depth4pt width0pt}%
    \normalbaselines\rm}

% this macro defines a twelvepoint font

\def\twelvepoint{
    \def\rm{\fam0\twelverm}
    \textfont0=\twelverm \scriptfont0=\ninerm \scriptscriptfont0=\sevenrm
    \textfont1=\twelvei \scriptfont1=\ninei \scriptscriptfont1=\seveni
    \textfont2=\twelvesy \scriptfont2=\ninesy \scriptscriptfont2=\sevensy
    \textfont3=\tenex \scriptfont3=\tenex \scriptscriptfont3=\tenex

    \textfont\itfam=\twelveit \def\it{\fam\itfam\twelveit}
    \textfont\slfam=\twelvesl \def\sl{\fam\slfam\twelvesl}
    \textfont\bffam=\twelvebf \def\bf{\fam\bffam\twelvebf}
    \textfont\ttfam=\twelvett \def\tt{\fam\ttfam\twelvett}
    \textfont\Bbbfam=\twelveBbb \def\Bbb{\fam\Bbbfam\twelveBbb}
    \scriptfont\Bbbfam=\nineBbb \scriptscriptfont\Bbbfam=\sevenBbb
    \def\sc{\twelvesc}
    \tt \ttglue=.5em plus.25em minus.15em
    \normalbaselineskip=14pt

    \setbox\strutbox=\hbox{\vrule height10pt depth4pt width0pt}%
    \normalbaselines\rm}

%\newfam\Bbbfam \def\Bbb{\fam\Bbbfam\twelveBbb} \textfont\Bbbfam=\twelveBbb
%

%%%%%%%%%%%%%%
%%
% this macro defines a elevenpoint font
%
%
\def\elevenpoint{
   \def\rm{\fam0\elevenrm}
   \textfont0=\elevenrm \scriptfont0=\eightrm \scriptscriptfont0=\sevenrm
   \textfont1=\eleveni \scriptfont1=\eighti \scriptscriptfont1=\seveni
   \textfont2=\elevensy \scriptfont2=\eightsy \scriptscriptfont2=\sevensy
   \textfont3=\tenex \scriptfont3=\tenex \scriptscriptfont3=\tenex

   \textfont\itfam=\elevenit \def\it{\fam\itfam\elevenit}
   \textfont\slfam=\elevensl \def\sl{\fam\slfam\elevensl}
   \textfont\bffam=\elevenbf \def\bf{\fam\bffam\elevenbf}
   \textfont\ttfam=\eleventt \def\tt{\fam\ttfam\eleventt}
   \textfont\Bbbfam=\elevenBbb \def\Bbb{\fam\Bbbfam\elevenBbb}
   \scriptfont\Bbbfam=\eightBbb \scriptscriptfont\Bbbfam=\sevenBbb
   \tt \ttglue=.5em plus.25em minus.15em
   \normalbaselineskip=12pt

   \setbox\strutbox=\hbox{\vrule height10pt depth4pt width0pt}%
   \normalbaselines\rm }

%%%%%%%%%%%%%%%%%%%%%%%%%%%%%%
% Invoke the default font size
\twelvepoint

% Set some default dimensions
%
\abovedisplayskip 14pt plus 3pt minus 10pt%
\belowdisplayskip 14pt plus 3pt minus 10pt%
\abovedisplayshortskip 0pt plus 3pt%
\belowdisplayshortskip 8pt plus 3pt minus 5pt%
\parskip 3pt plus 1.5pt
\hsize=6.5in
\vsize=8.9in

\fi			% end of \ifx\UsualIsLoaded
%%%%%%%%%%%%%%%end of usual.tex%%%%%%%%%%%%%%%

\noindent {\twelvepoint {\bf 0. Introduction}}
\tenpoint

\bp
A finite rank hermitian  super
vector bundle  ${\cal E}= {\cal E}^+ \oplus {\cal E}^-$ has an associated
determinant bundle
 $\Det ({\cal E})\equiv \left(\Det\ {\cal E}^+ \right)^*\otimes \Det\ {\cal E}^-$. A
connection
 $\nabla^{\cal E}$ on ${\cal E}$ with curvature $\Omega^{\cal E}$
 induces a connection
$\nabla^{\Det\ {\cal E}} $  on the determinant bundle, with curvature 
$\Omega^{\Det\ {\cal E}} = -\str(\Omega^{\cal E})$ equal to
minus the first Chern form  
on the
 original bundle
${\cal E}$.

 \sp In this paper, we investigate whether this property carries over
 to  infinite rank bundles of physical interest. 
The immediate problem is that $\str(\Omega^{\cal E})$ involves a
 divergent sum.
The paper breaks the problem down into two parts: (i) constructing the
 determinant bundle associated to an infinite rank superbundle,
 following [BF], [Q1]; 
(ii)  defining the first Chern form of the superbundle, and relating it to 
 the curvature
on the determinant bundle. 

\sp As background, Quillen [Q1] 
  constructed  the determinant bundle with a natural metric associated to
a family of Cauchy-Riemann operators on a Riemann surface, and
computed its curvature.
Later, Bismut and Freed [BF]
equipped the determinant bundle associated to a family of Dirac-type
operators with a
connection compatible with this Quillen
metric, and computed the curvature in terms
of local
invariants of the underlying spin manifold.  Freed [F] considered
characteristic forms on loop groups, overcoming divergence problems
via an {\it ad hoc} summation technique.  In [AP], [MRT], more natural
(but less tractable) heat kernel and zeta function regularization
techniques were used to renormalize divergent expressions.

\sp In \S\S1-3,
we formalize the construction of
Quillen-Bismut-Freed determinant bundles in terms of
determinant bundles associated
to    ``half  weighted super vector bundles.''  We first restrict ourselves
 to a class of   super
vector bundles ${\cal E}\equiv  {\cal E}^+\oplus {\cal E}^-$, where ${\cal
E}^\pm$
are vector bundles with fibers 
modeled on Sobolev spaces $H^{s^\pm}(M, E^\pm)$
 of sections of some finite rank
hermitian vector bundles $E^\pm$ over a 
closed Riemannian manifold $M$. A {\it half-weighted vector bundle}
 is such a hermitian superbundle together with a field/family
 $L\equiv 
\left[\matrix{ 0&L^-\cr
L^+ &0 \cr}\right]$ of odd self-adjoint  operators  locally given by
elliptic
 operators  acting on  smooth sections of $E\equiv  E^+\oplus E^-$. 
This local
 characterization makes sense globally if  the transition maps are
themselves zero order, grading preserving
 elliptic   operators   on $M$. To a half
weighted super vector bundle
$({\cal E}, L)$
we associate a determinant bundle $ \Det ({\cal E},L)$, the
Quillen determinant bundle of the family $  L^+ $.

\sp Given a half-weighted vector bundle $({\cal E}, L)$, we have a
family
 $Q\equiv  L^2= L^-L^+\oplus L^+L^- $  of
positive, self-adjoint, locally elliptic  operators acting fiberwise on
${\cal E}$.
As in [Pa], we call  $({\cal E},
Q)$   a
{\it weighted vector bundle}.  The {\it weight} $Q$ can be viewed as 
metric data on the infinite dimensional vector bundle ${\cal E}$, and
 the existence of
$L$ allows us to view ${\cal E}$ as a spinor bundle with 
Clifford multiplication given by the {\it half weight} $L$.

\sp
Starting in \S4, we construct regularized first Chern forms.
Using $Q$, we 
define $Q$-{\it weighted traces} $\tr^Q$ and $Q$-{\it  weighted
supertraces} $\str^Q$, which are linear functionals  on  sections of
$PDO({\cal E})$, the bundle of operators which are locally
 given by classical pseudo-differential operators on the fibers of ${\cal E}$.
We define the {\it weighted first Chern form}
of a superconnection $\nabla^{\cal E}$
on $({\cal E}, Q)$  as the
$Q$-weighted supertrace
 $\str^Q(  \Omega^{\cal E}  )$
of the
 curvature of the connection, provided $\Omega^{\cal E} $ is a
two-form with values in pseudo-differential operators on the fibers of
${\cal E}$.  

\sp
Our main results (\S6, Theorems 3, 5) show
 that  the curvature of the Bismut-Freed connection 
on the determinant bundle associated
to a
half-weighted superbundle  with connection
differs from   (minus) 
the weighted first Chern form on the superbundle
by  a linear combination of
 Wodzicki residues. This   obstruction to the finite dimensional formula
arises from the nonvanishing of $[\nabla^{\cal E},\str^Q]$, a
feature of the infinite dimensional weighting procedure.
We    express this obstruction in two ways:

\sp --  via zeta function regularization, using weighted
supertraces and evaluating the obstruction
 $[\nabla^{\cal E}, \tr^Q]$ in
terms of a
 Wodzicki residue (Theorem  3);

\sp -- via heat kernel regularization, using a
one-parameter family of Bismut
connections [B],
thus avoiding  weighted supertraces  (Theorem 5).

\sp\noindent We also show (Corollary 6) that the weighted first Chern
form is more local than the curvature of the Bismut-Freed connection in
a certain technical sense.  In the proof of the Corollary, we see that
the curvature of the superbundle is a multiplication operator and
therefore 
not trace-class.  Thus regularization procedures are necessary to
define the first Chern form.

\sp  In \S7, we  illustrate the main results with a string theory/Teichm\"uller
theory example.  Here the action of $H^{s+1}$ diffeomorphisms of a
closed surface $\Lambda$ on the manifold ${\cal A}(\Lambda) $ of
smooth almost complex structures on $\Lambda$ gives rise to a family
$\alpha_J: H^{s+1}(T\Lambda)\to H^s(T_1^1 \Lambda), J \in {\cal
A}(\Lambda) $ of elliptic operators.  Setting ${\cal E}^+\equiv T
{\cal A}^s( \Lambda) \mid_{{\cal A}(\Lambda)}$ and ${\cal E}^-\equiv
{\cal A}(\Lambda) \times H^s(T_1^1\Lambda)$, we can view $\left( {\cal
E} \equiv {\cal E}^{+ } \oplus {\cal E}^{- } , L\equiv
\left[\matrix{ 0&\alpha^*\cr
\alpha &0 \cr}\right] \right)$ as a half-weighted
superbundle. We identify the Bismut-Freeed superconnection with  classical
connections in Teichm\"uller theory.

\sp In Appendix A, we collect some superconnection calculations.  
 In Appendix B, as suggested by
the different proofs of Theorems 3 and 5, we relate Wodzicki residues
to the trace forms of [JLO].

\medskip
 $\bullet${\bf Notation:} Let $E$ be a
finite rank hermitian or Riemannian vector bundle over a Riemannian
manifold $M$. The natural $L^2$ inner product on the smooth sections of
$E$ is defined by
$$\langle \sigma, \tau\rangle\equiv  \int_M  \langle \sigma(x), \tau(x) \rangle_x d
\mu(x),$$
 where $\mu$ is the volume measure on $M$, and
$\langle \cdot, \cdot \rangle_x$ the
inner product on the fiber of $E$ above $x$.

\sp We
  denote by $PDO(M, E)$ the algebra of classical pseudo-differential
operators (PDOs) 
acting on smooth sections of $E$, by $Ell(M, E)$ the multiplicative
subset of
 elliptic  PDOs, by $Ell^{s.a} (M, E)$ the subset of self-adjoint elliptic
PDOs and  by
   $Ell^+(M, E)$   the subset of positive elliptic PDOs. Adding the
subscript $\ord >0$ to these sets
 restricts to operators of strictly positive order.  Adding the
superscript $*$
restricts to injective operators.

\sp
In the following we take $s>{\dim M\over 2}$. Recall that for  $s>{\dim
M\over 2}$, we have
$H^{k+s}(M, E)\subset C^k(M, E)$ for any $k\in \N$ where $ H^t(M,E)$ (resp.
$C^k(M, E)$) denotes
  the space of $H^t$ (resp. $C^k$) sections of the bundle $E$.

 \bp
\noindent {\twelvepoint {\bf 1. A class of vector
bundles}}

\tenpoint
 \np  We say that a
  Hilbert space $H$ lies in the class ${\cal C}{\cal H} $   if
there
is   a closed smooth Riemannian manifold $M$,
  a finite rank hermitian/Riemannian  vector bundle $E$ over $M$,
 and  $s>{\dim \ M \over 2}$ such that
  $H=H^s(M,E)$.
For example, for $G$ be a Lie group and ${\rm Lie}(G)$ its Lie algebra,
 the Lie
 algebra
$H^s(M, {\rm Lie} (G))$ of the Hilbert current group $H^s(M, G)$
 lies in ${\cal C} {\cal H}$.

\sp
 Let
 ${\cal C} {\cal E}$ be the class of Riemannian
Hilbert vector bundles   ${\cal E}\to X$
over a
  (possibly infinite dimensional) manifold $X$
with    fibers    modeled
on a separable Hilbert space $H=H^s(M,E)$ in 
${\cal C}{\cal H} $    and  with  
transition maps   in $PDO(M, E)$.
Note that these PDOs
have coefficients only in
some Sobolev class.  However, the PDOs in the examples below
are locally given by
multiplication operators, and are as tractable as PDOs with smooth
coefficients.  

 \sp
 ${\cal C}  {\cal X}   $  denotes the class
of   infinite dimensional manifolds $X$
 with
tangent bundle $TX$   in ${\cal C} {\cal E} $. Since the transition maps
are bounded,
 they correspond to operators of order zero.
Moreover,  the transition maps are invertible, so they
in fact lie
 in $Ell(M, E)$.

\np We now give 
examples of   manifolds in  ${\cal C}
{\cal X}$ and vector bundles in
 ${\cal C}{\cal E}$.

 \np $\bullet$ {\bf Examples:}\sp
\noindent i)  Finite rank vector bundles  lie in 
${\cal C}{\cal E}$. To see
this, we take
  as base manifold   a point $\{\ast\}$, and as the bundle $E$
  the trivial bundle
 $\{\ast\}\times \R^d$ (or $\{\ast\}\times \C^d$ if the bundle is complex).
The transition functions
belong to $Ell (\{\ast\} , E)= Gl_d (\R)$ (or $Gl_d (\C)$).
 We say that $M$ is reduced to a point.

\sp\noindent  ii) If $G$ is a Lie group and $s>{\dim\ M \over 2}$,  the
current group
 $H^s(M, G)$ is a Hilbert Lie group having
a left invariant atlas $\phi_\gamma(u)(x)\equiv   \exp_{\gamma(x)}
( u(x)),$ for $x\in M, \ \gamma \in H^s(M, G)$,
where $\exp_{\gamma(x)}$ is the exponential coordinate chart at $\gamma(x)$
induced by a
left invariant Riemannian metric  on $G$. The transition functions
are given by multiplication operators, which indeed are PDOs.

 \sp \noindent
iii) Let $M\equiv \Lambda$  be a closed,  oriented,  Riemannian 
surface of genus $p>1$, and let ${\cal
A}^s (\Lambda), s>1,$ be the space of almost complex structures on
$\Lambda$ of Sobolev class $H^s$, i.e.
$${\cal A}^s (\Lambda) = \{ J\in H^s (T^1_1\Lambda), J^2_x = -\hbox{Id}_x,
J_x \ \hbox{preserves\ orientation\ of}\ T_x\Lambda \ \hbox{for}
\ x\in\Lambda\}.$$ 
   ${\cal A}^s (\Lambda)$ is a smooth
Hilbert manifold  with tangent space at $J \in {\cal A}^s
(\Lambda)$  given by [T]
$$T_J {\cal A}^s(\Lambda) = \{ H\in
H^s (T^1_1 \Lambda), \ HJ  + J H=0\}.$$
(The set of smooth almost complex structures ${\cal A}(\Lambda)=
\bigcap_{s>1} {\cal A}^s(\Lambda)$ is only a Frechet manifold.)
We determine the transition maps. The charts are given pointwise by
the matrix exponential map
$\exp_JH(x)\equiv  \exp_{J(x)}H(x)$. Hence the transition maps as maps on
$H^s(T^1_1\Lambda)$ are multiplication operators,
  so  they
 are    PDOs of
order zero. Thus $T{\cal A}^s(\Lambda)$ is in ${\cal
C} {\cal E}$ with fibers modeled on $H^s(\Lambda,E)$ 
where $E\equiv   T_1^1\Lambda$.
In the string theory example in Appendix B,  we consider the subbundle
given by restricting $T{\cal A}^s(\Lambda)$ to the manifold
 ${\cal A}(\Lambda)$:
 $${\cal E}^-\equiv T{\cal A}^s(\Lambda)\mid_{ {\cal A} (\Lambda) }.
\eqno (1.1)$$  
 ${\cal E}^-$ has an almost
complex structure
defined fiberwise by $${\cal J}^-(J) (H)\equiv  J \cdot H,$$
where $\cdot$ denotes pointwise matrix multiplication. Notice that if $J$
is smooth  and $H$ of class $H^s$, then $JH$ is of class $H^s$.
${\cal J}^-$ induces a splitting
$${\cal E}^-\equiv  {\cal E}^{-^{1,0}}\oplus  {\cal E}^{-^{ 0,1}},
$$
where the fibers of the subbundles   above   $J\in {\cal
A}(\Lambda)$ are
$${\cal E}_J^{-^{1,0}}\equiv  {\rm Ker}({\cal J}^-_J-i), \quad {\cal E}_J^{-^{ 0,
1}}\equiv  {\rm Ker}({\cal J}^-_J+i).$$
Because the almost complex structure is defined pointwise by
$(J\cdot H)(x)= J(x)
 H(x)$ for $x\in \Lambda$ and hence defines a PDO, the transition functions
of
these subbundles are also given by PDOs. Thus
${\cal E}_J^{-^{1,0}},\ {\cal E}_J^{-^{0,1}}$
   lie  in  ${\cal C} {\cal E}$.

 \bp
\noindent  iv) In Appendix B, we also consider the trivial bundle
$$ {\cal E}^+\equiv     {\cal A} (\Lambda) \times H^{s+1}(T  \Lambda
),\eqno(1.2)$$ 
which clearly lies in ${\cal C} {\cal E}$.
${\cal E}^+$ has a natural almost complex structure ${\cal J}^+$ defined
fiberwise by the almost complex structure on the tangent space to
$\Lambda$:
$${\cal J}^+ (J) u\equiv  Ju.$$
With respect to the complex structure $J$,  $T\Lambda$
splits into
  $T\Lambda= T^{1,0} \Lambda\oplus T^{0,1}\Lambda$,
 with $T^{1,0} \Lambda\equiv 
{\rm Ker}(J-i),
\ T^{0,1} \Lambda\equiv  {\rm Ker}(J+i)$. 
 ${\cal E}^+$ therefore splits into 
subbundles,
  $$ {\cal E}^+= {\cal E}^{+^{1,0}}\oplus {\cal E}^{-^{1,0}}, $$
whose fibers above
 $J\in {\cal A}(\Lambda)$ are
 $${\cal E}_J^{+^{1,0}}\equiv  \Ker({\cal J}^+_J-i)=H^{s+1} (\Ker (J-i)),
 \quad {\cal E}_J^{+^{ 0, 1}}\equiv  \Ker({\cal J}^+_J+i)=H^{s+1} (\Ker
(J+i)).$$

  \bigskip

\noindent {\twelvepoint {\bf 2. Weighted vector bundles and half-weighted super
vector bundles }}

\tenpoint
     \np A {\it weighted Hilbert space} is a pair $(H, Q)$ with $H$ in
${\cal C}{\cal H} $
and
 $Q\in Ell_{\ord >0}^{+ }  (M, E)$.
   \np $\bullet$ {\bf  Bundles of elliptic operators:} Let ${\cal E}$
 be a vector bundle in  ${\cal C}{\cal E} $  over
 a manifold $X$  with ${\cal E}$
modeled on a separable Hilbert space
$H$.
 For $x\in X$, let $PDO({\cal E}_x)$ be
 the set of operators
  $A_x $ acting densely on the fiber ${\cal E}_x$ above $x$
  such that    for any local trivialization
$\phi: {\cal E}|_{U_x}  \to U_x\times  H$ near $x$, the operator
$ \phi^\sharp A(x)\equiv  \phi(x)A_x \phi(x)^{-1} $  lies in $PDO(M,
E)$.
Here $\phi(x): {\cal E}_x
\to H$ is the
isomorphism induced by the trivialization.
Similarly, let $Ell({\cal E}_x)$ be the set of    operators
$A_x $ acting densely on
$T_xX$
 such that     for any local trivialization $\phi: {\cal E}|_{U_x}  \to U_x
\times  H$
near $x$, 
the operator $\phi^\sharp A    $ lies in $Ell(M, E)$. 
From this point on we will omit the subscript $x$.

  \sp  These definitions are  independent of the
  choice of local chart. Indeed, since transition functions are given by
operators in $PDO(M, E)$,
the condition
$\phi^\sharp A \in PDO(M, E) $
 is independent of the choice  of $\phi $.
Since the principal symbol is multiplicative and since
ellipticity is characterized by invertibility of the principal symbol,
the condition  $\phi^\sharp A(x) \in Ell(M, E) $
 is also independent of the choice of   $\phi $. Notice that the order of
$\phi^\sharp A$ is independent of the choice of local chart, so 
 we can speak of the order of $A$.

 \sp This gives rise to  bundles $PDO( {\cal E}),\ Ell({\cal E})$
with fiber at  $x$ given
respectively by
  $PDO(
 {\cal E}_x),
\ Ell({\cal E}_x)$.  In particular, a section of the second bundle is
a family of elliptic operators parametrized by the base.
When
 ${\cal E}$ is a bundle of finite rank, we can view it  as before
as a bundle of sections over
a manifold reduced to a point. Then $PDO({\cal E}) =
\Hom({\cal E},{\cal E})$ and
 $Ell ({\cal E})= Gl({\cal E})$. 

\sp
 \np $\bullet$ {\bf Weighted  bundles:}   A local
  section $Q$
of $Ell ({\cal E})$, with ${\cal E}$  modeled on some $H^s(M, E)$,
 is {\it
positive self-adjoint} if for all
 $x$ in the support of $Q$, and in any local chart $(U, \phi)$
around $x$, the operator
 $  \phi^\sharp Q(x)   $  
  lies in $Ell^+(M, E)$. 
 A {\it weighted   bundle}  is a
pair $({\cal E},
Q)$ with ${\cal E}$ in
 ${\cal C}{\cal E}$ and $Q$ a section of positive
self-adjoint operators of constant order in
$  Ell({\cal E})$.
 A {\it weighted manifold} $(X, Q)$ is a manifold in  ${\cal C}
{\cal X}$ such that
 $(TX, Q)$ is a weighted vector bundle.
The operator  $\phi^\sharp Q $ is by definition
a {\it weight} on the model space $H$ of $X$.

    \np
$\bullet$ {\bf Examples:}  We return to examples i)-iv).

 \sp i) A Riemannian structure on a finite rank vector bundle is a
weight, since it
yields a family of positive definite linear transformations,
which are
positive self-adjoint elliptic operators acting on sections of a vector
bundle over a manifold
reduced to a point.  

 \sp ii)  For the
current groups
$H^s(M, G)$,
   let
$Q_0\equiv 
\Delta \otimes 1_{{\rm Lie}(G)}$ be the Laplace-Beltrami 
operator on $M$ with values in
the Lie algebra ${\rm Lie} (G)$ of the group $G$, for $\Delta$ the 
Laplace-Beltrami operator acting
on complex valued functions on $M$.  For $\gamma \in H^s(M, G)$,
setting 
$Q(\gamma)\equiv  L_\gamma^{-1} Q_0 L_\gamma$, where $L_\gamma$ is left
multiplication by $\gamma$,
yields a weighted manifold
 $(H^s(M, G), Q)$. 

\sp
 iii) and (iv)  We consider the bundles ${\cal E}^\pm$ defined above.
 For $J\in {\cal
A}  (\Lambda)$, let $\alpha_J : H^{s +1}(T\Lambda)\rightarrow H^{s } (T^1_1
\Lambda)$ be the operator defined by the Lie derivative of $J$:
$$\alpha_J u = {d\over{dt}} \big( f^*_{u,t }J\big),\eqno(2.1)$$
where $f_{u,t}$ is the flow of the
 vector field $u$.
$\alpha_J$ is a first order
elliptic operator with range   $T_J {\cal A}^s (\Lambda) = \{ H\in H^{s }
(T^1_1 \Lambda), HJ + JH = 0\}$ [T]. Its adjoint  $\alpha_J^*  $ is
defined 
with respect to
  the hermitian products:
$$\eqalignno{\langle
u,v\rangle_{J }^+ &\equiv  \int_\Lambda d\mu_J (x) \ 
%u\cdot_{g_J} v
\langle u,v\rangle_{g_J}, &( 2.2^+)\cr
\langle
H,K\rangle_{J }^- &\equiv  
\int_\Lambda d\mu_J (x) \ 
%\hbox{tr} \ (HK^*)
\langle H, K\rangle_{g_J}.&( 2.2^-)\cr}$$
Here $g_J$ is the unique metric of constant curvature $-1$ among the
conformal class of metrics for which $J$ is orthogonal [T], and
$d\mu_J$ is the associated volume form.  Note that $\langle H,
K\rangle = {\rm tr}(HK^*)$, where $K^*$ is the hermitian adjoint of
the matrix representing the $(1,1)$ tensor $K$ 
with respect to $g_J.$  Since
%in isothermal coordinates.
$\alpha^*_J \alpha_J$ and $\alpha_J \alpha^*_J$ are elliptic,
 the families 
$$Q^+ \equiv 
 \{ Q_J^+\equiv  \alpha^*_J
\alpha_J, J\in {\cal A} (\Lambda)\},\ 
Q^-\equiv \{ Q^-_J\equiv  \alpha_J \alpha^*_J, J \in {\cal
A}(\Lambda)\},$$
yield  weighted bundles
$ ({\cal E}^+, Q^+),\ 
({\cal E}^-, Q^-),$ respectively.  Thus
 we get  a weighted super vector bundle:
$$\left( {\cal E}= {\cal E}^+ \oplus {\cal E}^-, Q\equiv  Q^+ \oplus
Q^-\right).\eqno(2.3)$$

 \np
 $\bullet${\bf Half-weighted super vector bundles:}  For a   super vector
bundle $ {\cal E} $
in  ${\cal C}{\cal E}$ with fibers modeled on some
$H^{s^+}(E^+)\oplus H^{s^-} (E^-), s^\pm >{\dim\ M\over 2}$, via local
charts we can write
 a local section
$L$ of 
 $Ell({\cal E})$    in  matrix form
$L=\left[\matrix{ L_{++} & L_{+-}\cr
L_{-+} & L_{--} \cr}\right]$. Provided the transition maps are even,
it makes sense to consider the  class of odd operators,
i.e.~those which locally have   only off-diagonal terms.
 We define 
   a  {\it   half-weighted superbundle } to be  a pair  $({\cal E},    L )$,
where ${\cal E}$ is a superbundle in
$ {\cal C} {\cal E}$ with even transition maps and $L$ is 
a    section of 
odd self-adjoint  operators in $Ell ({\cal E})$
of non-zero order.

\np To a half-weighted superbundle $({\cal E},   L )$ we can associate
 a weighted superbundle
$({\cal E}, Q\equiv   L^2)$. Since $L$ is odd, we can write
  $L\equiv   \left[\matrix{ 0 & L^{- }\equiv \left( L^+\right)^*\cr
L^{+ } & 0 \cr}\right]$, so the weight $Q$ can be written as
$$ Q  = Q^+  \oplus Q^-   \equiv    L^{-}   L^{+}\oplus L^{ +} {L^{ -}} 
=  {( L^{  +})}^*   L^{ +}\oplus {( L^{  -})}^* {L^{  -}}. $$

$\bullet$ {\bf Examples:}
$\left({\cal E}\equiv {\cal E}^+\oplus {\cal E}^- ,
 L\equiv \left\{\left[\matrix{ 0 &   \alpha_J^*\cr
\alpha_J & 0 \cr}\right], J \in {\cal A}
(\Lambda)\right\} \right)$, with ${\cal E}^\pm$ as in (1.1), (1.2),
 and $\alpha_J$ as in (2.1), is a half-weighted super vector
bundle. If
 $\alpha_J$ stablizes the fiber
${\cal E}_J^{1,0}$  for each $J\in {\cal A}(\Lambda)$, 
we can build a complex half-weighted
bundle
 $({\cal E}^{1,0}, L^{1,0}\equiv \{L_J^{1,0}, J \in {\cal A}(\Lambda)\})$,
where $L_J^{1,0} = \alpha_J|_{{\cal E}^{1,0}_J}.$

\sp As shown in Appendix B,   $L_J^{1,0}$  is  a
Cauchy-Riemann operator,
the historically first case of examples
  provided by
 spinor bundles on   even dimensional manifolds [Q1], [BF].
 Let  $\pi : Z\rightarrow B$ be a smooth fibration of even
dimensional spin 
manifolds  $\{M_b, b\in B\}$, and let ${\cal E}\rightarrow B$ be
an infinite dimensional super 
vector bundle with fiber $H^s(M_b,E_b)$ for
a smooth family $\{ E_b,
b\in B\}$ of Clifford bundles on $M_b$. The
Dirac operators $D_b = D_b^+ \oplus D_b^-$ 
act on $H^s(M_b,E_b)$ as elliptic operators. 
For
$$  L_b \equiv 
\left[ \matrix { 0 &D_b^- =  ({D_b^+}^*)\cr
D_b^+ &0 \cr}\right], \eqno(2.4)$$
 $({\cal E},L)$
is a half-weighted superbundle.

\np $\bullet$ {\bf  From group actions to half-weighted superbundles: }
Half-weighted superbundles also arise from group actions.
Let $\G$ and $\P$ be two infinite dimensional Hilbert manifolds modeled
respectively on $H^{s^+} (M, E^+)$ and $H^{s^-} (M, E^-)$, where $E =
E^+\oplus E^-$ is a superbundle over $M$, such that:
 \item{a)} $\G$ has  a smooth group   multiplication on the right:
$
R_{\gamma_0}:\G\rightarrow \G,\ 
\gamma\mapsto\gamma \gamma_0$,
for $ \gamma_0\in \G$.
\item{b)} $\G$ acts on ${\cal P}$ on the right by
$
 \Theta  : \G\times \P\rightarrow \P,\ 
(\gamma, p)\rightarrow p.\gamma,$
inducing
 a smooth map
$
\theta_p :\G\rightarrow \P,\ 
\gamma\rightarrow p.\gamma$, for
$ p\in \P.$
\item{c)} The differential
 $\alpha_p \equiv  d\theta_p : T_e\G\rightarrow T_p$
 is  elliptic, with  order independent of
$p$.

\sp \noi Let  ${\cal E}^+ \equiv   B\times \Lie (\G)$, where  $B$ is a
submanifold of $\P$, $\Lie({\cal G}) = T_e{\cal G}$,
and ${\cal E}^-= T{\cal P}|_B.$
Then
 $$\left({\cal E} = {\cal E}^+ \oplus {\cal E}^-, {\cal L} = \left\{ L_b\equiv 
\left[\matrix{ 0&
 \alpha_b^*\cr
\alpha_b &0\cr}\right], b\in B\right\}\right) $$
is a half-weighted superbundle.
 
\sp
$\bullet$ {\bf Example:}
In  the notation of Examples   (iii) above, let
${\cal G}\equiv {\hbox{Diff}}^{s+1}_0 (\Lambda)$
be the group of  isotopies (i.e.~diffeomorphisms homotopic to the
identity) of
$\Lambda$ of Sobolev class $H^{s+1}$.
Although ${\cal G}$ is
not a Lie group, it is a Hilbert manifold modeled on $H^{s+1} (T\Lambda)$
with a smooth multiplication on the right.  ${\cal G}$ acts on ${\cal A}^s (\Lambda)$
 (which
we recall is modeled on $H^s(T^1_1)$) by pullback, and this action satisfies
 a) and b) above (see   [T]).  Since 
$\alpha_J$ in (2.1)
 is elliptic, the family
$   L\equiv  \left\{\left[\matrix{ 0&
 \alpha_J^*\cr
\alpha_J &0\cr}\right], J\in {\cal A}\right\} $ yields a half-weighted
structure on
the bundle
 ${\cal E}  $ in (2.3).
  \vskip 1cm

\noi{\twelvepoint
{\bf 3. From a half-weighted  super vector bundle to the determinant
bundle}}
\tenpoint

\np    Let $({\cal E},     L )$   be a half-weighted super vector 
bundle over  a manifold $B$, which as above determines the  weighted
superbundle 
$({\cal E},  Q)$ with
$Q = L^2$.
From $({\cal E},    L )$  we construct the
 {\it determinant bundle $\Det ({\cal E}, L) = \Det({\cal E})$}, 
following [BF], [BGV],
[Q1].

As before, set ${\cal E}= {\cal E}^+\oplus {\cal E}^-$  where
${\cal E}^\pm$ has fibers modeled on $H^{s^\pm}(M, E^\pm)$, and  write
a section $L$ of $Ell({\cal E})$ consisting
of odd, self-adjoint operators in the form
$\left[\matrix{ 0 & L^-= ({L^+})^* \cr
L^+  &0\cr }\right].$
 Let $m$ be the  order of $L_b, b\in B,$ and set
$s_+ = s_- +m$. This 
yields a family of Fredholm operators $L_b^+:  H^{s^+}(M, E^+)\to  H^{s^-
}(M, E^-)$.
As Quillen shows,
there is a  line bundle, the determinant bundle $\Det({\cal E})$
over $B$, with fiber $\Det({\cal E})_b \simeq
(\Lambda^{\rm top}\Ker \ L^+_b)^* \otimes \Lambda^{\rm top}
\Coker\ L^+_b$, where $\Lambda^{\rm top}$ denoes the top exterior
power.  $\Det({\cal E})$ has a canonical section Det $L^+$ given by
$(e_1\wedge\ldots\wedge e_n)\otimes (f_1\wedge\ldots\wedge f_m)$,
where $\{e_i\}$, resp.~$\{f_j\}$ are orthonormal bases of the
eigenvalues of $L_b^-L_b^+$, resp.~$L_b^+L_b^-$, lying below some
$a\in {\bf R}$ not in the spectrum of either operator.

  \sp 
$\bullet$ {\bf A family of  connections on the determinant bundle:} Fix
$\varepsilon>0$.  At any point $b\in
B$ where $L_b$ is injective, the $\varepsilon$-cutoff determinant of
the self-adjoint elliptic operator $Q^+_b = L^-_b L^+_b$ is defined by
$$\hbox{det}_\varepsilon Q^+_b \equiv  \hbox{exp} \left[
-\int^{\infty}_\varepsilon {1\over t}{\hbox{tr} (e^{-tQ^+_b}}) \
dt\right].$$ 
These yield a one-parameter family of Quillen metrics
$\{\Vert\cdot\Vert_{Q,\varepsilon}, \varepsilon > 0\}$ on $\Det({\cal E})$
defined by
$$\Vert\hbox{Det}\ L^+_b\Vert_{Q,\varepsilon} \equiv 
\sqrt{\hbox{det}_\varepsilon Q^+_b},$$
if $L_b^+$ is invertible.
At other points, we replace $Q_b^+$ in the previous equation
 by $Q_{b,a}^+$,
the restriction of
$Q_b^+$ to the eigenspaces above $a$.
Given a connection $\nabla^{\cal E}$ on ${\cal E}$, as in
[BF] we can define a one-parameter family $\{
\nabla^{\Det({\cal E}) ,\varepsilon}, \varepsilon > 0\}$ of
connections on $\Det({\cal E})$
compatible with the metrics
$\{\Vert\cdot\Vert_{Q,\varepsilon}, \varepsilon > 0\}$ by
$$\eqalign{ (\Det\ L^+_b)^{-1} \nabla^{\Det({\cal E}),\varepsilon}
\Det\ L^+_b
&\equiv  \hbox{tr} \Big((L^+_b)^{-1} \nabla^{ \Hom({\cal E}) } L^+_b 
e^{-\varepsilon Q^+_b}\Big)\cr
&= {1\over 2}\left(d\log\ \det {}_\e Q_b^+ + 
\hbox{str}\left( (L_b)^{-1} \nabla^{ \Hom({\cal E}) } L_b 
e^{-\varepsilon Q_b}\right)\right).\cr}\eqno(3.1)$$
Here $\nabla^{ \Hom({\cal E})} $ denotes  the connection on Hom$({\cal E}^+,
{\cal
E}^-)$ induced by $\nabla^{\cal E}$, and str denotes the supertrace,
defined by
$$\str(A) = \hbox{str}\left[\matrix {A^+ & X\cr
Y & A^- \cr} \right]\equiv  \tr(A^+)-\tr(A^-).\eqno (3.2)$$
This definition is motivated by the corresponding formula
for the natural connection (5.3)
on the determinant line bundle for a finite
rank superbundle.
In (3.1) and from now on, we assume that $L_b^+$ is injective.
In general, our
formulas can be modified via the cutoff operator $Q_{b,a}^+$.

\np $\bullet${\bf Renormalized limits:} 
Following [BGV, Ch.~9], from the family  of 
   connections $\{\nabla^{\Det({\cal E}),
 \e}, \e>0\}$,  we build  a renormalized
connection   by taking a    renormalized limit   as
$\e\to 0$. 
  More precisely, 
%   $f: \R^+\setminus\{0\}\to \C$ have the asymptotic expansion as $\e\to 0$
%   $$f(\e)=  \sum_{i=-K}^J \alpha_j \e^{j\over m} +  \sum_{q=0}^N \beta_q
%  \e^q \log\ \e+
%  {\rm O}(\e^{J+1\over m})$$
%   for some $J\in \N$, $m\in \N\setminus\{0\}$,$ N\in \N$ and $\alpha_j, \beta_q \in
%  \C$. For $\mu \in \R$,
%   the {\it $\mu$-renormalized limit} at $0$ is defined to be
%  $$\Lim^\mu_{\e \to 0} f(\e)\equiv  \alpha_0-\mu \beta_0.\eqno(3.3)  $$
for $(m, n)\in (\N\setminus\{0\})\times \N$,
$\alpha\in \R$,
  let   ${\cal F}_{m, n,\alpha}$ be the set
of functions   $f: \R^{+}\setminus \{0\}\to \C$ such that there exist
$a_j,   b_j, c_j, j\in \C$
with
  $$f(\e)\sim\sum_{j=0 }^\infty a_j \e^{\l_j }
+ \sum_{j=0,\l_j\in \Z }^\infty  b_j\e^{\l_j }\hbox{log}\ \e + \sum_{j=0
}^\infty
c_j\e^j,  $$
as $\e\to0$,
where $\l_j\equiv {j-\alpha-n \over m} $. In other words, 
 for $J \in \N$ and
 $ K_J\equiv [\alpha]+m J+n\in \N $, we have
$$f(\e)=\sum_{j=0  }^{K_J}   a_j \e^{\l_j }
+ \sum_{j=0  ,\l_j\in \Z}^{K_J}  b_j\e^{\l_j } \hbox{log}\ \e +
\sum_{j=0}^J
c_j \e^j
 + {\rm o}(\e^J)  $$
(cf.~(4.3)).
(If $\alpha\in \Z$, there is a redundancy 
since constant terms can arise in the first
and last sum.)
We call such a function {\it  renormalizable},
 and for  $f\in {\cal F}_m\equiv
 \bigcup_{n\in \N, \alpha \in \R} {\cal F}_{m,n,
\alpha}$,  we define the {\it  $\mu-$renormalized limit  } of $f$ at
zero   by:
$$ \Lim^\mu_{\e \to 0} f(\e)
= a_{\alpha+n} + c_0- \mu b_{\alpha+ n}, \eqno(3.3) $$
where  we set $a_{\alpha+n}=0$   and $b_{\alpha+n}=0$ if
$\alpha +n \not \in \N$.  Thus $\Lim^\mu_{\e \to 0} f(\e)$ is the
constant term in $f$'s asymptotic expansion minus $\mu$ times the
coefficient of $\log\ \e.$  ([BGV] only consider the case $\mu=\gamma,$
the Euler constant.)
If there is no logarithmic divergence, then $\Lim\equiv  \Lim^\mu
\equiv \Lim^\mu_{\e\to 0}$ is
independent of $\mu$.

 \np $\bullet$ {\bf A renormalized connection on the determinant bundle:}
As in [BF], [BGV], 
we set
$$\eqalign{ (\Det\ L^+_b)^{-1}
\nabla^{ \Det ,\mu}
\Det\ L^+_b
&\equiv   \Lim_{\e \to 0}^\mu\hbox{tr} \Big((L^+_b)^{-1} (\nabla^{ \Hom({\cal E}) }
L^+_b )
e^{-\varepsilon L^-_b L^+_b}\Big)\cr
&= {1\over 2}\left( d\log\ \det 
{}_\mu Q_b^+ + \Lim_{\e \to 0}^\mu 
\hbox{str}\left( (L_b)^{-1} (\nabla^{ \Hom( {\cal E}) } L_b )
e^{-\varepsilon Q_b}\right)\right),\cr}$$
where $$\hbox{det}_\mu Q^+_b \equiv  \hbox{exp} \left(
-\Lim_{\e \to 0}^\mu\int^{\infty}_\varepsilon {1\over t}{\hbox{tr} 
(e^{-tQ^+_b})} \
dt\right).  $$  
The renormalized connection $\nabla^{ \Det({\cal E}) ,\mu} =
\nabla^{ \Det,\mu}$ is compatible
with the renormalized
Quillen metric  given by
 $$\Vert\Det\ L^+_b\Vert_{Q,\mu } \equiv 
\sqrt{\hbox{det}_\mu Q^+_b}. $$
The curvature of $\nabla^{\Det,\mu}$ is denoted by
$\Omega^{\Det, \mu}$.

 \bp\noi
{\twelvepoint{\bf 4. First Chern forms on  weighted vector bundles}}
\tenpoint

\sp \noi The first Chern form on a finite rank hermitian bundle
with  connection is the trace of the curvature.  In
infinite
rank, one cannot expect  curvature to be trace-class in general,
so we
need to regularize (or renormalize) the trace.   We will use extra
data of the weights
of \S2 to define
weighted
traces in two steps: (i) defining a one-parameter
family of
weighted traces; (ii)  taking a renormalized limit.

\sp \noi    Let 
$({\cal E},    Q )$ be a weighted vector bundle in ${\cal C}
{\cal H}$
with fibers modeled on $H^s(M, E)$, and
let $ A $ be a section of  $PDO({\cal E}) $.
$Q$ is positive elliptic with strictly positive order, so for
$\e>0$,
$e^{-\e Q}$ is infinitely smoothing when seen in a local chart.  Thus
$Ae^{-\e Q}$ is trace-class when considered in a local chart
as a trace-class operator acting on $L^2 (M, E)$. 
We remark   that  a trace-class operator for the
$L^2$ inner product $\langle \cdot, \cdot \rangle$
 can be considered equally well as
 a trace-class operator with respect to the $H^s$ scalar
product
$$\langle \sigma, \rho\rangle^s\equiv 
\langle (Q+ 1)^{s/ \ord\ Q}\sigma,(Q+ 1)^{s/\ord\ Q} \rho\rangle. 
$$

 \np
$\bullet${\bf A family of weighted pseudotraces:} We define  a
 one-parameter family of   {\it $Q$-pseudotraces} of $A $
by
$$\hbox{tr}^{Q}_\varepsilon (A )\equiv  \hbox{tr}  (A
e^{-\varepsilon{  Q} }),\eqno(4.1) $$
for $\varepsilon >0$.
 Again this definition  should be understood in a local chart, but it
is easily
independent of the choice of chart, since for an
invertible operator $C$, we have
$\tr( CAC^{-1}e^{-\e CQC^{-1}})= \tr(C Ae^{-\e Q}C^{-1})=\tr(Ae^{-\e
Q})$.

  \np\noi 
We emphasize that pseudotraces  are not traces in the usual sense.
First,
$\tr_\e^{  Q} [A_1,A_2]\neq 0$ in general.
Moreover, unlike the  finite dimensional
case, if
$\{Q_t, t\in
\R\}$, is a one-parameter family of weights and
$  \{A_t, t\in \R\}$ is 
a one-parameter family of   PDOs, then for
fixed $\e>0$,
$${d\over dt}\biggl|_{_{t=0}}\tr_
\e^{Q_t}(A_t)\neq \tr_\e^{Q_0} ( \dot A ) -  \e
\tr^{Q_0}_\e(A_0 \dot Q )
$$
(where $\dot T\equiv  {d\over dt}|_{t=0} T_t$), as one would
expect
from a   formal differentiation of  (4.1), since  
in general neither $\dot Q$ nor $A_0$
commutes with $Q_0$.    (If either
$\dot Q$  or $A_0$ commutes with $Q_0$, this equation holds
by 
[G, \S1.9].)
These  obstructions can be analyzed more carefully using the  renormalized
pseudotraces in the
 next paragraph (see [CDMP], [Pa]).

\np  By (4.1), the one-parameter family of connections
on the determinant bundle  given by
(3.1) is
$$\eqalign{ (\hbox{Det}\ L^+ )^{-1} \nabla^{ \Det\ {\cal E},   \varepsilon}
\hbox{Det} \ L^+ &=
\hbox{tr}_\varepsilon^{Q^+ }  \big( (L^+ )^{-1} \nabla^{ \Hom({\cal E})}
L^+\big)\cr
&= {1\over 2}\left( d\log\ \det {}_\e Q^+ +
\hbox{str}_\varepsilon^{Q  } \big( L^{-1} \nabla^{ \Hom({\cal E})}
L\big)\right).\cr}\eqno(4.2)$$
As before $Q^\pm\equiv  L^\mp L^\pm$ and $Q\equiv  Q^+\oplus Q^-$.

\np
$\bullet${\bf  Renormalized  pseudotraces:} 
From the classical theory of  heat expansions
[G], [K], [Le, (3.18)], 
for   a positive elliptic operator $Q $ of positive integer order  and
 a PDO $A$ acting on sections of a
vector bundle over a closed manifold $M$,  the map $\e \to
\tr(A e^{-\e Q})$ lies in the class
${\cal F}_q$ of \S3, where $q = \ord(Q).$
More precisely,   there  exist
$N = N(\dim\ M) \in\N^+$, $a = a(\ord(A))\in\R$
and  $\alpha_j(Q,A),\ \beta_j(Q, A),\ \gamma_j(Q, A)\in\C$
%for  $j \in \N$,    
such that
  $$\tr(A e^{-\e Q}) \sim \sum_{j=0 }^\infty   \alpha_j(Q, A)
\e^{j-a-N\over q}
+ \sum_{j=0,  {j-a-N\over q}\in \Z}^\infty  \beta_j(Q, A)\e^{j-a-N\over
q}
\hbox{log}\  \e+ \sum_{j=0}^\infty \gamma_j(Q, A)  \e^j
 \eqno (4.3) $$ 
as $\e\to 0.$
(As before, 
constant terms can arise as both $\gamma_0$ and
$\alpha_{a+N}$ if $\ord(A)\in\Z$.)
%  For a PDO 
%  $A$  and
%    a positive    elliptic
%  operator $Q$ of positive integer order  acting on
%  sections of a  vector bundle over a compact manifold
%   $M$, we have the following asymptotic expansion for
%  $\tr(Ae^{-\e Q})$ as $\e \to 0$ (see [G], [K]):
%  $$\tr(Ae^{-\e Q})=  \sum_{j=- K }^J \alpha_j(A,Q) \e^{j/ \ord\ Q} +
%   \sum_{q=0}^N \beta_q(A,Q)\e^q \log\ \e  +{\rm O}
%  (\e^{(J+1)/ \ord \ Q})\eqno (4.3)$$
%  for some $K,N\in \N$,   $\alpha_j(A,Q), \beta_q(A,Q) \in \C$.
  We  define  the {\it $Q$-renormalized trace} of 
  $A$ as the $\mu$-renormalized limit  
  of the map $\e \mapsto \tr(Ae^{-\e Q})$, as in (3.3):
  $$\tr^{Q, \mu} (A)\equiv  \Lim_{\e \to 0}^\mu \tr_\e^{Q} (A)=
   \alpha_{a+N}(Q,A)+\gamma_0(Q,A)-\mu \cdot
  \beta_{a+N}(Q,A).\eqno (4.4) $$
  If we set
  $\mu=0$, resp.~$\mu= \gamma$, the Euler constant, we get a
   heat kernel renormalized trace, resp.~a 
  zeta function renormalized trace,  and these two are related via a
  Mellin  transform.
  In the following, we will usually consider the
  case $\mu=0$, and write $\tr^Q$ for $\tr^{Q, 0}$.
  The results can  easily be extended to the general case of
$\tr^{Q, \mu}$. 

\sp
 $\bullet$ {\bf Renormalized supertraces:}  (4.4) extends to
$Q$-renormalized supertraces in the obvious way for 
$Q=Q^+\oplus Q^-$:
$$\str^{Q, \mu}(A)=  \tr^{Q^+, \mu} (A^+) - \tr^{Q^-, \mu} (A^-), $$ 
for $A$ as in (3.2) with
$A^\pm$ PDOs.
%If $A$ is an ordinary pseudo-differential operator, namely $A=A^+$ then
%$str^{Q, \mu} (A)= tr^{Q, \mu} (A^+)$.  
Renormalized  pseudotraces/supertraces
appear in the  geometry  of determinant
bundles [BF], where the connection  on the determinant bundle 
can be written as
$$\eqalign{ (\hbox{Det}\ L^+_b)^{-1} \nabla^{ \Det ,\mu}
\hbox{Det}\ L^+_b
&\equiv    \hbox{tr}^{Q^+, \mu} \Big((L^+_b)^{-1} \nabla^{ \Hom( {\cal E}) } L^+_b
\Big)\cr
&= {1\over 2}\left( d\log\ \det {}_\mu Q^+ +
 \hbox{str}^{Q , \mu}\left( (L_b)^{-1} \nabla^{ \Hom( {\cal E}) }
L_b \ \right)\right).\cr}
 \eqno(4.5)  $$
They also have been used 
(i) to define  minimality of infinite dimensional 
submanifolds of
manifolds of connections and metrics [AP], [MRT], and (ii)
in relation
to determinants of elliptic operators 
[KV 1],  for a   special class of operators  on
which they are
actually  traces. 

\np These renormalized traces are related to Wodzicki residues, as we
briefly recall; see
 [KV1, 2], [Pa] for  more details.
Let $({\cal E},    Q )$ in ${\cal C}
{\cal H}$ be a  weighted vector
 bundle
with fibers modeled on $H^s(M, E)$, and
  let $ A $ be a section of  $PDO({\cal E}) $.
Since $Q$ is positive elliptic with strictly positive order, for
any $z\in \C$ with
${\rm Re}(z)>{\dim \ M / \ord\ Q}$, the
 operator  $(Q+ P_Q)^{-z} $ is trace-class on $L^2 (M, E)$
in any local
chart. Here $P_Q$ is the orthogonal projection onto the kernel of $Q$.
 Similarly,
  for $ {\rm Re}(z)> {
 (\dim\ M+ \ord \ A) / \ord\ Q}$, 
 $A(Q+ P_Q)^{-z} $ is trace-class.
For such $z$, we may define
    $$\widetilde{  \str}_z^{Q}  (A)\equiv  \hbox{str}  (A
  {(Q  +P_Q)}^{-z}).$$
By the Mellin transform, we have
$$\beta_0(Q,A)={\rm res}_{z=0}(\widetilde {\str}_z^{Q}(A)) $$
in the notation of (4.3); in particular, $\beta_0(A) = \beta_0(Q,A)$
is independent of $Q$.  It follows that
$$\str^{Q, \mu}(A)= \lim_{z \to 0} \left( \widetilde{ \str}_z^{Q}(A)-
z^{-1}\res_{z=0}(\widetilde{ \str}_z^{Q}A) )+ (\gamma-\mu )
\res_{z=0}(\widetilde 
{\tr}_z^{Q}(A)\right). $$
 Renormalized pseudotraces thus arise as the finite part  of a divergent
expression. The infinite part is built from
the Wodzicki residue [W] $\res(A)$:
 $$\res(A)\equiv  (\ord\ Q) \cdot  \res_{z=0}
(\widetilde{ \tr}_z^{Q}
 (A) ),\eqno(4.6)$$ 
which
 defines a trace on the algebra of pseudo-differential operators
[W], [K]. In summary:
$$\str^{Q, \mu}(A)= \lim_{z \to 0} \left( \widetilde {\str}_z^{Q}(A)-
{1\over z \cdot \ord\ Q }\res ( A )\right)+ { \gamma-\mu \over \ord\ Q} \res
  (A ).   $$

\bp \noi
  We can now define  $Q$-weighted first Chern forms
on a weighted vector
 bundle.

 \sp \noi
{\bf Definition:} {\it  Let $({\cal E}, Q)$ be a weighted
hermitian   (super) vector bundle over $B$ 
with  connection $\nabla^{\cal E} $
and curvature
$\Omega^{\cal E}$.
Assume that for any  $X, Y  \in \Gamma({\cal E})$,   
$ \Omega^{\cal E}(X, Y) \in\Gamma(PDO({\cal E}))$.  Define
\sp {i)} the one-parameter family of  {\rm $ Q$-weighted first Chern
forms} by
$$r_{1 }^{Q,\e} (X,Y)\equiv \str_\e^{Q}\left(\Omega^{\cal E}(X,Y) \right),\quad \e
>0,  \eqno (4.7)$$
 \sp {ii)} the one-parameter family of
   {\rm $Q$-renormalized first Chern forms}
$$R_1^{Q, \mu}(X,Y)\equiv   \str^{Q, \mu}\left(\Omega^{\cal E}(X,Y)\right), \quad
\mu \in
\R. \eqno (4.8)$$}

  \bp \noi
{\twelvepoint{\bf 5. 
The curvature on the associated determinant bundle in finite
dimensions} }
\tenpoint

  \np
Let ${\cal E}$ be a finite dimensional bundle with connection
$\nabla^{\cal E}$, and let
 $\alpha$ be a   $\Hom( {\cal E},{\cal E})$-valued form.
Writing
 $\nabla^{\cal E}= d+ \theta$ in a local trivialization, we have:
   $$ d\ \tr(\alpha) =   \tr([d,
\alpha] ) 
=\tr([d,\alpha])+ \tr([\theta, \alpha])
=\tr([\nabla^{\cal E}, \alpha]),\eqno (5.1 )$$
since the trace term  $\tr([\theta, \alpha])$ vanishes.
The final expression is of course independent of the
 choice of local trivialization.
Thus the trace of a covariantly constant form is closed.
In particular, since the curvature $\Omega^{\cal E}$
is covariantly constant by the Bianchi
identity, the first Chern form
$r_1^{\cal E}\equiv  \tr(\Omega^{\cal E})$ 
is also closed.
This form is a representative of the first Chern
class in de Rham cohomology.

 \sp
This generalizes to  supertraces on superbundles:
$$ d\ \str(\alpha) = \str([\nabla^{\cal E}, \alpha]), \eqno (5.2)$$
where $[\cdot, \cdot]$ is now a supercommutator and $\nabla^{\cal E}$ a
superconnection on the superbundle ${\cal E}.$
The first Chern form 
$r_1^{\cal E}\equiv  \str(\Omega^{\cal E})$
is therefore also closed.

 \sp We recall the relation between the first
 Chern form of a superbundle and the curvature of
the associated determinant bundle.
Let ${\cal E}^\pm$ be
hermitian vector bundles   with connections $\nabla^{{\cal
 E}^\pm}$ over a manifold $B$. 
$\nabla^{{\cal
 E}^\pm}$ induce a connection $\nabla^{\cal E}$ on ${\cal E} =
 {\cal E}^+ \oplus {\cal E}^-.$
%Let us write
%$\nabla^{{\cal E}^+}= d +\theta^+$ and $\nabla^{{\cal E}^-}= d+ \theta^-$
%in a local  trivialization. Then the homomorphism bundle
The bundle 
$\Hom({\cal E}^+, {\cal E}^-)\simeq ({\cal E}^+)^*\otimes {\cal E}^-$
has the natural connection
$\nabla^{\Hom(\cal E)}\equiv  {(\nabla^{{\cal E}^+} )}^*\otimes 1+ 1\otimes
\nabla^{{\cal E}^-} $, given by $\nabla^{\Hom({\cal E})}L^+ =
[\nabla^{\cal E},L^+]$ for $L^+\in\Gamma(\Hom({\cal E}^+, {\cal
E}^-))$ (cf.~Appendix A).
% Locally we have $\nabla^{\cal E}= d+ \theta\equiv  d+ \theta^--
%\theta^+$.
Assuming for convenience that ${\cal E}^\pm$  have the  same rank,
the determinant bundle
$\Det({\cal E})\equiv  (\Lambda^{\rm top} 
{\cal E}^+)^* \otimes
\Lambda^{\rm top} {\cal E}^-$
has the
hermitian metric
$$\Vert \Det \ L^+\Vert\equiv  \sqrt{\det ((L^+)^*L^+)}$$
for $L^+\in\Gamma(\Hom({\cal E}^+, {\cal E}^-))  $ and $\Det\ L^+ $
the corresponding section of
$\Det({\cal E}^+, {\cal E}^-)  $. 
 $\nabla^{\cal E}$ induces a connection $\nabla^{\Det\ {\cal
E}}$  on $\Det({\cal E})  $ compatible with this metric,
defined at points where $L^+$ is injective by:
$$\eqalign{ (\Det\ L^+)^{-1} \nabla^{\Det\ {\cal E}} \Det\ L^+&\equiv 
  \hbox{tr}(\left(L^+\right)^{-1}[\nabla^{\cal E}, L^+]).\cr
&= {1\over 2} \left(  d\log\ \det Q^+ +
  {\rm str}\left(L ^{-1}[\nabla^{\cal E},
L]\right)\right),\cr}\eqno (5.3)$$
where $L = L^+\oplus (L^+)^*$, $Q^+ = (L^+)^*L^+$
(cf.~(4.3), (4.5)).  The following lemma
is well known.

  \bigskip
\noindent\hang {\bf Lemma 1:} {\it  The curvature   $\Omega^{\Det\ {\cal E}}$
of the connection $\nabla^{\Det\ {\cal E}}$ on the
 determinant bundle $\Det({\cal E})$ associated to the 
connection
$\nabla^{\cal E}$ on the
superbundle
${\cal E}= {\cal E}^+\oplus {\cal E}^-$ satisfies}
$$\Omega^{\Det\ {\cal E}} = -\str(\Omega^{\cal E})
= {\rm ch}(\nabla^{\cal E})_{(2)},\eqno(5.4) $$
{\it where $\Omega^{\cal E}$ is the curvature of $\nabla^{\cal E}$,
and ${\rm ch}(\nabla^{\cal E})_{(2)}$ is the degree two component of the Chern
character $\str(\exp[-\Omega^{\cal E}])$ of the connection:
i.e.~the curvature of the determinant line bundle is minus the first
Chern form of the superbundle.}
%and for any invertible section
%$L^+\in\Gamma(Hom({\cal E}^+, {\cal E}^-))$, we have}
%$$(\Det\ L^+)^{-1}\Omega^{\Det\ {\cal E}}  \Det\ L^+=
%\str(\Omega^{\cal E} ). \eqno (5.7)$$
\bigskip

\noindent {\bf Proof:}  For later purposes, we give a basis free
proof.  
Pick $M,N\in T_bB$, where $L^+$ is injective
at $b$.  Extend $L^+$ near $b$ so that $[\nabla^{\cal E},L^+]_b=0.$
By (5.3), we have
$$\eqalign{
%[(\Det\ L^+)^{-1}\Omega^{\Det\ {\cal E}}  \Det\ L^+](M, N)&=
\Omega^{\Det\ {\cal E}}(M,N) &=
{1\over 2}   d\left({\rm str}(L ^{-1}[\nabla^{\cal E}, L])\right)
(M, N)\cr
&= {1\over 2} \str\left(\left[\nabla^{\cal E}, L^{-1}[\nabla^{\cal E},
L]\right]\right)(M,N).\cr}$$
Using the Cartan formula $d\alpha(M,N) = M(\alpha(N))-N(\alpha(M))
-\alpha([M,N])$, we get
$$\eqalignno{
%[(\Det\ L^+)^{-1}\Omega^{\Det\ {\cal E}}  \Det\ L^+](M, N) 
\Omega^{\Det\ {\cal E}}(M,N)
&= {1\over 2}\biggl(-{\str}\left(\left[L ^{-1}[\nabla_M^{\cal E}, L] ,L
^{-1}[\nabla_N^{\cal E}, L]
\right] \right) + {\str}
\left(L ^{-1}[\Omega^{\cal E}, L]\right)(M, N)\biggr)\cr
&= {1\over 2} {\str}\left(L ^{-1}[\Omega^{\cal E}, L]\right)(M, N) &(5.5)\cr
&= -  {\str}( \Omega^{\cal E}) (M, N),\cr}$$
where we have used  
$\str(A^{-1}[B, A])= -2\cdot\str(B)$ for $A$ odd, $B$ even.
The second equality in (5.4) is standard.\hfill$\bullet$

       \bp
\noindent{\twelvepoint   
 {\bf  6. The curvature on the determinant bundle in 
infinite dimensions}}
\tenpoint

 \sp The main goal of this paper is to see how 
(5.4)   extends to the infinite
dimensional
setting.
More precisely, the Quillen-Bismut-Freed theory of
determinant bundles constructs a determinant bundle with connections
(4.3), (4.5), for
certain half-weighted superbundles, with the curvature 
of (4.5) computed in [BF].  Via weighted traces, we have
constructed weighted and renormalized first Chern forms of
such superbundles, and it is natural to ask if
(5.4) continues to hold.  

\sp The proof of (5.4)
uses the facts $\tr([A,B])=0$ and 
$d\ \str = \str([\nabla,\cdot])$, both of which fail for weighted traces.
Thus we cannot expect (5.4)
 to hold in infinite dimensions.
 Indeed we  will show by two methods that
 (5.4) holds up to an
obstruction given by
 Wodzicki residue terms defined in (4.6). 
The two methods lead to different expressions for
these obstructions which
 seem difficult to identify  directly.

\sp The first  zeta function regularization approach   uses weighted
traces to express the supertrace of a commutator and the obstruction to
$d\ \str = \str([\nabla,\cdot])$
in terms of
 Wodzicki residues.  The appearance of Wodzicki residues is natural,
since they are defined via zeta function regularization.
 The second  heat kernel regularization approach  uses a one-parameter
family of
 superconnections introduced by Bismut [B] to avoid 
weighted traces,  and  closely follows the methods used in [BF], [BGV] to
compute the curvature on the
 determinant bundle for families of Dirac operators.

\np $\bullet$ {\bf  First approach using weighted traces} 

\sp
  Let $ {\cal E}= {\cal E}^+ \oplus {\cal E}^-$ in ${\cal C}{\cal E}$ 
be a superbundle 
with  connection $\nabla^{\cal E}$
 with even
transition maps acting on a model space $H^s(M, E)$, $s>{\dim\ M\over 2}$.
Let $Q$ be a weight on ${\cal E}$.
 The following lemma expresses the obstruction to $d\ \str =
\str([\nabla,\cdot])$ 
as a Wodzicki residue.

\bigskip
\noindent\hang {\bf Lemma 2:} [CDMP] {\it 
  Let $({\cal E},Q)$ be a weighted vector
bundle with connection
$\nabla$ over a manifold $B$, 
and let $\alpha$, $\beta$  be $PDO({\cal E} )$-valued
one-forms on $B$.
 For $\mu \in \R$,}
\sp 1)     $\displaystyle{\str^{Q, \mu}[\alpha, \beta]= 
-{1\over \ord\ Q} \res ([\log\ Q,
\alpha] \beta).}$\hfill (6.1)\break
\sp 2)   {\it if $[\nabla, \log\ Q]$  and $ [\nabla, \alpha]$ are
$PDO({\cal E})$-valued one-forms, then}
$$   d\left(\str^{Q, \mu}(\alpha)\right) = \str^{Q,
\mu}([\nabla,\alpha])-{1\over \ord\ Q} \res(\alpha \cdot [\nabla,
\log\ Q]). \eqno(6.2)$$

  \sp  
For completeness, we 
outline the proof of (6.2) for traces, which easily extends to
supertraces, and refer the reader to [CDMP] for (6.1).
As before, $\tr^Q$ denotes the
renormalized trace $  \tr^{Q, \mu}$
at $\mu=0$; the results extend to 
 $\mu\neq 0$.

 \sp  One first shows that for one-parameter families of operators
$A_t\in PDO(M, E),$  $Q_t\in Ell^+_{\ord >0}(M, E)$
of constant order, 
we have
$$   {d\over dt}\biggl|_{_{t=0}}
 \left(\tr^{Q_t}(A_t)\right) = \tr^{Q_0}\left({d\over
dt}\biggl|_{_{t=0}}A_t\right)-{1\over \ord\ Q_0}
 \res\left(A_0  {d\over dt}\biggl|_{_{t=0}} \log\ Q_t\right).
$$
This uses the fundamental property of the canonical trace of
Kontsevich-Vishik [KV 2].
Similarly, in a fixed local trivialization of ${\cal E}$,
we have
$$   d\ \tr^Q(\alpha) = \tr^Q(d\alpha )-{1\over \ord\ Q} \res(\alpha \cdot d
\log\ Q  ). \eqno(6.3)$$
Let  $\nabla= d+ \theta$ in the local trivialization.  Since 
$[\nabla, \alpha]=d \alpha + [\theta, \alpha]  \in
\Gamma( PDO({\cal E}))$, 
and  since $d\alpha$, 
the
 differential of a PDO, also lies in $\Gamma( PDO({\cal E}))$, 
it follows that
$ [\theta, \alpha]$ lies in $PDO(M, E)$ pointwise.
 Using again the fundamental
property of the canonical trace, one shows
$$\tr^Q[\theta, \alpha]= -{1\over \ord\ Q} \res ( [\log\ Q, \theta]
\alpha).\eqno(6.4)$$
Combining (6.3) and (6.4) gives
$$ \eqalign{   d\ \tr^Q(\alpha) &= \tr^Q(d\alpha )-{1\over \ord\ Q} \res(\alpha
\cdot d
\log\ Q  )\cr
&= \tr^Q([\nabla, \alpha] )- \tr^Q ([\theta, \alpha])-{1\over \ord\ Q}
\res(\alpha \cdot d
\log\ Q  )\cr
&= \tr^Q([\nabla, \alpha] )+ {1\over \ord\ Q} \res([\log\ Q, \theta]\alpha
)-{1\over \ord\ Q} \res(\alpha
\cdot d
\log\ Q  )\cr
&= \tr^Q([\nabla, \alpha] ) -{1\over \ord\ Q} \res(\alpha
  [\nabla,
\log\ Q] ).\cr}$$
${}$\hfill$\bullet$

 \bp  The
 residue term in (6.2) is the source of the infinite dimensional 
obstruction to identifying
 the first Chern form of a superbundle with (minus) the curvature of the
 determinant bundle:

 \bp
\noindent\hang {\bf Theorem 3:  } {\it   
%The curvature   $\Omega^{Det,
%\mu}  $  on the
% determinant bundle   associated to the half-weighted  superbundle
Let $\left({\cal E}= {\cal E}^+\oplus {\cal E}^-, L = L^+\oplus
L^-\right)$  be a half-weighted  super vector bundle
with connection $\nabla^{\cal E}$ over a manifold $B$.
The curvature   $\Omega^{\Det,
\mu}  $  of the associated 
 determinant bundle differs from
   the $Q$-weighted first Chern form $ R_1^{Q, \mu}$
of (4.8)  on the weighted superbundle  $\left({\cal E}, Q=L^2\right)$
by a Wodzicki residue.
  More precisely,  
%if $L_b^+$ is injective, 
for $M,N\in T_bB$ we
have}
$$ 
%    \left( \Det\ L^+ \right)^{-1} \Omega^{\Det, \mu}(M,N)  \Det\ L^+  =
 \Omega^{\Det, \mu}(M,N) =
  -R_1^{Q, \mu}(M, N)  + {\cal R}^{Q, \nabla^{\cal E}} (M, N),$$
{\it with  }
$$\eqalign{ 2{\cal R}^{Q, \nabla^{\cal E}} (M, N)&=   {1\over \ord\ Q
} \res
\biggl(
 \left[ \log\ Q ,  L^{-1} [\nabla_M^{\cal E}, L ]\right]
  L  ^{-1}[ \nabla_N^{\cal  E},
 L ]
\cr
 &\qquad -    L  ^{-1} [\nabla_N^{\cal E}, L ]
\cdot [\nabla_M^{\cal E},
\log\ Q ]\cr
 &\qquad +  L^{-1} [\nabla_M^{\cal E}, L ]
\cdot [\nabla_N^{\cal E}, \log\ Q ]\biggr). \cr}
$$

 \sp \noindent {\bf Proof:} 
 We follow the proof of Lemma 1, replacing traces by
renormalized supertraces and keeping track of obstructions due to (6.1) and
(6.2) via Wodzicki residues. We obtain
 $$\eqalign{
%2 (\Det\ L^+)^{-1}\Omega^{\Det}(M, N) \Det\ L^+&=
2\Omega^{\Det}(M, N)  &=
d \left(\hbox{str}^{Q }( L^{-1} [\nabla^{\cal E}, L])\right)(M, N)\cr
 &=
 -\hbox{str}^{Q }\left(
 \left[L^{-1} [\nabla_M^{\cal E}, L],   L^{-1}[ \nabla_N^{\cal  E},
L] \right]\right) \cr
&  \qquad + \str^{Q }\left(L^{-1} [\Omega^{\cal E}(M, N), L]\right)\cr
&\qquad -{1\over \ord\ Q} \res(   L^{-1} [\nabla_N^{\cal E}, L] \cdot [\nabla_M^{\cal
E}, \log\ Q])\cr
 &\qquad +{1\over \ord\ Q} \res(   L^{-1} [\nabla_M^{\cal E}, L]  \cdot
[\nabla_N^{\cal E}, \log\ Q]),  \cr}$$
using (6.2) and calculating as in (5.5).  Thus
$$\eqalign{2 \Omega^{\Det}(M, N) 
 &=-2 \ \str^{Q}\left( \Omega^{\cal E}(M, N) \right)\cr
  &\qquad +{1\over \ord\ Q} \res\left(
 \left[\log\ Q,  L^{-1} [\nabla_M^{\cal E}, L]\right]
L^{-1}[ \nabla_N^{\cal  E} , L ]\right) \cr
 &\qquad -{1\over \ord\ Q } \res( L^{-1} [\nabla_N^{\cal E}, L ] \cdot
[\nabla_M^{\cal E}, \log\ Q ])\cr
 &\qquad +{1\over \ord\ Q } \res(    L^{-1} [\nabla_M^{\cal E}, L]
\cdot [\nabla_N^{\cal E}, \log\ Q ]), \cr}$$  
using (6.1).\hfill$\bullet$

\bp $\bullet$ {\bf   The heat kernel approach:}   Here we deform the
weight $Q = Q_0\equiv L^2$ to 
a one-parameter
family $Q_0+ Q_{1, \e}$, $\e>0$ via a deformation of the
superconnection $\nabla^{\cal E}$
into  a   family $\nabla^{L}_\e$ of Bismut
superconnections.
We need a preliminary formula.

 \np $\bullet$ {\bf Volterra series:} [BGV, (2.5)]
Let $Q=Q_0+Q_1$, where $Q_0$ is a
positive elliptic operator of
strictly positive order, and $Q_1$ is a PDO
 of order strictly less than that of $Q_0$.
We have
$$    e^{-\e(Q_0+ Q_1)}  = \sum_{k=0}^\infty
(-\e)^k\int_{\Delta^k}
%{{\sigma_i\geq 0\atop 
%\sum_{i=0}^k \sigma_i= 1}} 
e^{-\sigma_0
\e Q_0} Q_1 e^{-\sigma_1
\e Q_0} Q_1 \cdots Q_1 e^{-\sigma_k
\e Q_0} d\sigma_0d\sigma_1 \cdots d\sigma_k,  $$
where $\Delta^k = \{\sigma_0,\ldots,\sigma_k >0: \sum_{i=0}^k \sigma_i =1\}.$
 We can avoid convergence issues,
since we will only be using 
a finite number of terms.
In analogy with the notation in [JLO], we set
$$\langle A_0, A_1, \cdots, A_k\rangle_{  \e,k, Q_0}\equiv  \int_{\Delta^k}
 \str\left( A_0 e^{-\sigma_0
\e Q_0}  A_1 e^{-\sigma_1
\e Q_0} A_2 \cdots A_k e^{-\sigma_k
\e Q_0} \right)d\sigma_0d\sigma_1 \cdots d\sigma_k,   $$
for 
PDOs $A_0, \cdots, A_k$ acting
on sections of the model bundle $E$ of ${\cal E}$.
The supertrace is clearly finite for
$\e>0$.
 The Volterra formula implies
$$    \str\left(e^{-\e(Q_0+ Q_1)}\right)  =
 \sum_{k=0}^\infty (-\e)^k \langle 1, Q_1, \cdots, Q_1\rangle_{   \e,k,
Q_0}. $$

 \np $\bullet$ {\bf Bismut superconnections:} 
Starting from a
half-weighted  superbundle
 $({\cal  E}= {\cal E}^+ \oplus {\cal E}^-,    L)$ with a metric
superconnection
$\nabla^{\cal E}= \nabla^+ \oplus \nabla^- $,
%  compatible with $\gamma^0$,
we form the
one-parameter family
of superconnections
$$\nabla_\e^L\equiv  \nabla^{\cal E}  +\sqrt \e L,$$
for $\e >0$ [B]. 
For any one-parameter family of
superconnections $A_t$, we have the important {\it transgression
formula}: for an analytic function $f$, 
$${d\over dt} \str(f(A_t^2))= d\biggl(\str\left({d \over d
t} A_tf'
(A_t^2)\right)\biggr).\eqno(6.5)$$
The derivation of this formula in [BGV, Prop. 1.41] essentially relies
on the fact that $d\left(\str (\alpha)\right)= \str[A, \alpha]$ for
$\Hom({\cal E},{\cal E})$-valued one-forms
$\alpha$.  Thus in the proof below, we avoid the
obstructions to (5.1), (5.2), for the weighted supertraces
$\str^Q.$

\bigskip
\noindent \hang{\bf Proposition 4:}
{\it  Let $({\cal E}, L)$ be a half-weighted superbundle 
with  connection
$\nabla^{\cal E}$ over a manifold $B$.
For  $\e>0$  we have}
$$\Omega^{\Det\ {\cal E},\e}  =  \ch(\nabla_\e^L)_{[2]}=
-  r_1^{Q,\e}
+ {\e  } \langle  I, [\nabla^{\cal E}, L], [\nabla^{\cal E}, L] \rangle_{
\e, 2, Q_0},$$
{\it where  $Q_0=L^2$,  $r_1^{Q,\e}\equiv \str_\e^{Q_0}(\Omega^{\cal
 E})$
 is the
weighted  first Chern form of
$\nabla^{\cal E}$ as in (4.7), $\Omega^{\Det\ {\cal
E}, \e}$ is  the
curvature of $\nabla^{\Det\ {\cal E} , \e}$ as   in (3.1), and
$\ch(\nabla_\e^L)_{[2]} \equiv \str\left(\exp\left(
- (\nabla_\e^L)^2\right)\right)_{[2]}$
is the degree two component of the Chern character of
 $\nabla_\e^L$. }
\bigskip

\noindent {\bf Proof:} 
We first compute the degree  $k$ piece of the Chern character of
$\nabla_\e^L$ in two 
ways.
On the one hand,  since
$$\left( \nabla_\e^L\right)^2=
\e Q_0+ \sqrt \e [\nabla^{\cal E}, L]+
\left(\nabla^{\cal E}\right)^2= \e( Q_0+ Q_{1, \e}), $$
with 
$ Q_0\equiv  L^2,\  Q_{1, \e}\equiv  {1\over \sqrt \e }[\nabla^{ \cal E},
L]+{1\over \e}\left(\nabla^{\cal E} \right)^2$,
by the Volterra formula we have
$$\ch(\nabla_\e^L)_{[k ]}= \sum_{j=0}^\infty (-\e)^j\big[ \langle I, Q_{1, \e},
\cdots, Q_{1, \e}\rangle_{\e, j, Q_0}\big]_{[k ]}.\eqno(6.6)$$
As promised, we need only
consider
a finite number of terms   in this
sum. 

\sp
On the other hand,  we have
$$\eqalign{ \ch (\nabla_\e^L)_{[k]}&= \str \left[ \exp[-(\nabla_\e^L)^2]
\right]_{[k]}
= -\left[ \str  \int_\e^{\infty} \left( {d\over dt} ( \exp[-
(\nabla_t^L)^2])\right)\right]_{[k]} dt\cr
&= -\int_\e^{\infty} {d\over dt}\left[ \str \left(    \exp[-
(\nabla_t^L)^2] \right)\right]_{[k]} dt.\cr}$$
Applying the transgression formula (6.5)  to
$f(x) =   e^{-x}$, we get
$$\eqalign{\ch (\nabla_\e^L)_{[k]}
&=   \int_\e^{\infty}d \left[ \str \left(\left( {d\over dt}\nabla_t^L\right)
 \exp[-( \nabla_t^L)^2] \right)\right]_{[k-1]} dt\cr
&= {1\over 2} \left( d \left[ \int_\e^{\infty} \str \left({L\over \sqrt t}
\exp[-(\nabla_t^L)^2]\right) dt \right]_{[k-1]} \right)\cr
&= {1\over 2}  d \int_\e^{\infty} \sum_{j=0}^\infty {(-t)^j \over \sqrt t}
\big[ \langle L, Q_{1, t}, \cdots, Q_{1, t}\rangle_{ t, j ,
Q_0}\big]_{[k-1]}.\cr}\eqno(6.7)$$
Combining (6.6) and (6.7) yields
$$\sum_{j=0}^\infty (-\e)^j \big[\langle I, Q_{1, \e},
\cdots, Q_{1, \e}\rangle_{\e, j, Q_0} \big]_{[k]}= {1\over 2}
d \int_\e^{\infty} \sum_{j=0}^\infty {(-t)^j \over \sqrt t}
\big[ \langle L, Q_{1, t}, \cdots, Q_{1, t}\rangle_{ t, j ,
Q_0}\big]_{[k-1]}.\eqno(6.8)$$
Using Lemma 1 and recalling that $Q_0 = L^2$, we have
$$\eqalign{ \Omega^{\Det\ {\cal E}, \e}&= {1\over 2} d
\left(\str (L^{-1} [\nabla^{\cal E}, L] e^{-\e Q_0})\right)
= -{1\over 2}d \int_\e^\infty {d\over dt} \str( L^{-1} [\nabla^{\cal E}, L]
e^{-tQ_0}) dt \cr
 &=  {1\over 2} d\int_\e^\infty  \str( L  [\nabla^{\cal E}, L]
e^{-tQ_0}) dt 
 =  {1\over 2}  d \int_\e^\infty \langle L, [\nabla^{\cal E}, L] \rangle_{
t, 1, Q_0}dt\cr
 &=
{1\over 2}  d \left[\int_\e^\infty \sqrt t\langle L,  Q_{1, t}  \rangle_{
t, 1, Q_0}dt\right]_{[1]}
  =  -{1\over 2}  d \int_\e^\infty \big[\sum_{j=0}^\infty {(-t)^j\over \sqrt t}
\langle L,Q_{1, t}, \cdots, Q_{1, t}\rangle_{ t, j,
Q_0}\big]_{[1]}dt.\cr}\eqno (6.9)$$
 We used $[L,Q_0]=0$ in the second line, and the last equality
follows since the only term in the infinite sum of
 degree one is the integrand in the next to last integral.
By (6.7) and (6.9), we see that $\Omega^{\Det\ {\cal E}, \e} = 
-\ch (\nabla_\e^L)_{[2]}.$

\sp  Finally, by (6.8) and (6.9), we get
$$\eqalign{\Omega^{\Det\ {\cal E}, \e} &= 
-\langle I, \left(\nabla^{\cal E}\right)^2\rangle_{\e, 1, Q_0} + \e \langle
I, [\nabla^{\cal E}, L],
[\nabla^{\cal E}, L]\rangle_{\e, 2, Q_0}\cr
&= -\str(\Omega^{\cal E} e^{-\e Q_0}) +  \e \langle I, [\nabla^{\cal E}, L],
[\nabla^{\cal E}, L]\rangle_{\e, 2, Q_0},\cr}$$
which finishes the proof.\hfill $\bullet$

\bigskip
\noindent {\bf Remark:} In fact, (6.8) vanishes for $k$ odd, since the
integrand is the supertrace of an odd operator, and hence vanishes.
%   Setting  
%   $k=1$  in (6.8) yields
%  $$  -\sqrt \e\cdot \str \left([\nabla^{\cal E}, L] e^{-\e Q_0}\right) = {1\over
%  2} d \int_\e^\infty
%  t^{-{1\over 2}}\str\left(L e^{-tQ_0}\right) dt,\eqno (6.10)$$ 
%  a well known formula for the variation of the (cut-off)
%  $\eta$-invariant. In fact, for  
%  $$\eta_\e(L)\equiv   {1\over \sqrt \pi}\int_\e^\infty t^{-{1\over 2}}
%  \str( L e^{-t L^2} )dt,$$
%  (6.10) is equivalent to
%   $$d \eta_\e(L)= -{2\sqrt \e\over \sqrt \pi}  \str([\nabla^{\cal E}, L]
%  e^{-\e L^2}). $$
%  (This definition of the $\eta$-invariant as $\lim_{\e \to 0}
%  \eta_\e(L)$  coincides with the usual defintion
%  $$\lim_{z\to 0} {1\over \Gamma({z+1 \over 2})}\int_0^\infty t^{z-1\over 2}
%  \str( Le^{-tL^2}) $$ via a Mellin transform.)

\bigskip
By taking the renormalized limit in Proposition 4, we obtain:
\bigskip

 \np \noindent\hang
 {\bf Theorem 5:}  {\it For any $\mu\in \R$, the renormalized first
Chern  form $R_1^{Q, \mu}$
defined in (4.8) and the    curvature $\Omega^{\Det, \mu}$ of
the determinant line bundle are related by}
$$\Omega^{ \Det,\mu }  = \Lim^\mu_{\e \to 0}
\ch (\nabla_\e^L)_{[2]}=
   -R_1^{Q, \mu }   +\Lim^{\mu}\left(   \e \langle I, [\nabla^{\cal E}, L],
[\nabla^{\cal E}, L]\rangle_{\e, 2, Q_0}\right).$$

\bigskip

   \bp

\noindent $\bullet$ {\bf Bismut-Freed connections:}  For the
connection on the infinite dimensional bundle ${\cal E} = {\cal E}^+ \oplus
{\cal E}^-$   considered in
[BF], we can say more about the renormalized first Chern form. As in 
(2.4), we consider a fibration $\pi: Z\to B$ 
of  manifolds, with fiber an even dimensional manifold $M_b, b\in B$, 
and with finite
rank hermitian bundles $E^\pm$ with unitary connections
over $Z$.   
%  (The assumption in [BF] that the fibers are spin is not
%  necessary for our purposes.)
The Levi-Civita connection $\nabla^{\rm LC}$  for 
a given metric on $Z$ and the associated orthogonal horizontal 
splitting $T_xZ = T_bM_b \oplus H_x$, for $x\in \pi^{-1}(b),$ induces a 
connection $\nabla^F$ on $F$, the tangent bundle along the fibers of $\pi$ by
$$\nabla^F = P^{TM}\nabla^{\rm LC}, \eqno (6.10)$$
where $P^{TM}$ is the orthogonal projection to $F$.  We set 
${\cal E} = {\cal E}^+ \oplus {\cal E}^-$, a bundle over $B$ with
fiber 
$H^s(M_b, (F\otimes E^\pm)|_{M_b})$ for ${\cal E}^\pm.$  
The connections on $E^\pm$
induce a connection $\nabla^{F,E}$ on $F\otimes
(E^+\oplus E^-)$.  This lifts to a connection $\widetilde\nabla$ 
on ${\cal E}$
given by:
$$(\widetilde\nabla_Y h)(b)(x) = (\nabla_{\widetilde Y}^{F,E} h)(x),
\eqno(6.11)$$
where $\widetilde Y$ is the horizontal lift of $Y\in TB$ to $H_x$.
 This connection need not be
unitary, but it is shown in [BF] that adding the divergence of the
volume form at $x$ in base directions to the right hand side of (6.11)
produces a unitary superconnection $\widetilde\nabla^u,$  the Bismut-Freed
connection on ${\cal E}.$

We claim that the curvature form $\Omega^u$ for $\widetilde\nabla^u$ is an
endomorphism in the fibers: 
$$\Omega^u_b\in \Lambda^2(\Hom( F\otimes E|_{M_b},F\otimes E|_{M_b})),$$
for $E = E^+\oplus E^-$.
First, the local nature of the Bismut-Freed connection show
that if $\psi\in\Gamma({\cal E})$ has support in $U\times V$, where
$U$ is a  open set in $B$ and $V$ is an open set containing $x$ in
$\pi^{-1}(U)$, then $\widetilde\nabla^u\psi$ also has support in $U\times
V.$  This implies the same result for $\Omega^u(X,Y)$, a combination of
first and second covariant derivatives of $\widetilde\nabla^u.$
  Since $\Omega^u$ is tensorial,
after multiplying $\psi$ by bump functions with decreasing support in
base directions, we
can shrink $U\times V$ to the point $b$.  In other words,
$(\Omega^u(X,Y)\psi)(x)$ is determined by $\Omega^u(X,Y)_x$ and $\psi(x)$
alone.  Since $\Omega^u$ is linear, it must be an endomorphism in
the fibers.

This allows us both to
 compute the renormalized first Chern form $R_1^Q =
\str^Q(\Omega^u)$ and to relate $d\left(\str^Q(\Omega^u)\right)$ and
$\str^Q([\widetilde\nabla^u, \Omega^u]) =0.$  Here $Q=L^2$, with $L$ a 
first order differential operator as in (2.3).  Since $Q$ (which in 
[BF] is the square of the
Dirac operator in the fibers) is a second order differential operator,
we have the asymptotic expansion for the heat kernel $e^Q(\e,x,x)$:
$$e^Q(\e,x,x)  = {1\over (4\pi )^{(\dim \ M)/2}}\biggl(
\sum_{j = -\dim\ M}^J \alpha_j(x) \e^{j/2 } + {\rm
O}(\e^{(J+1)/2})\biggr).$$
$\alpha_j(x)\in \Hom(E_x,E_x)$ is locally computable
from the metric on $M_b$ at $x$ and the symbol of $Q$ at $x$.  
$\Omega^u$ is an endomorphism and $Q$ is a differential operator, so
$\str(\Omega^u e^{-\e Q})$  has an asymptotic expansion in $\e$ with no
logarithmic terms, and 
$\str^{Q}(\Omega^u)=\str^{Q,\mu}(\Omega^u)$ is consequently independent
of $\mu.$  In fact, by the  standard ``remarkable cancellations'' of local
index theory,
$\str(\Omega^ue^{-\e
Q})$  is ${\rm O}(1)$ as $\e\to 0.$ Thus we can replace $\Lim_{\e \to
0}^\mu $ in the definition of the renormalized first Chern form
$R_1^Q= R_1^{Q,\mu}$ by an ordinary limit:
$$\eqalign{R_1^Q &=\str^{Q}(\Omega^u) = \lim_{\e\to 0}
\str(\Omega^u e^{-\e Q}) 
=\lim_{\e\to 0}   \int_M \str(\Omega_x^u e^Q(\e,x,x))\cr
&= 
{1\over (4\pi )^{(\dim \ M)/2}}\int_M
\str^\mu(\Omega_x^u\alpha_0(x)),\cr}\eqno(6.12) $$
where we have used that $\Omega^u$ is a homomorphism
in the fibers.
Also,  
$$\eqalign{(4\pi)^{(\dim\ M)/2}
d\ \str^{Q}(\Omega^u) &= d\int_M \str(\Omega_x^u \alpha_0(x))
= \int_M \str([\widetilde\nabla, \Omega^u \alpha_0(x)])\cr
&= \int_M \str\bigl([\widetilde\nabla, \Omega^u]) \alpha_0 + \Omega^u [\widetilde
\nabla_x,\alpha_0]\bigr) = \int_M
\str(\Omega^u[\widetilde\nabla_x, \alpha_0]).\cr}$$
Thus the obstruction to $R_1^Q =\str^{Q}
(\Omega^u)$ being closed is given by the
integral of $\str(\Omega^u[\widetilde\nabla_x, a_0]).$  
The integrand is a local expression except in its dependence on 
$\widetilde\nabla_x.$  

\bigskip
\noindent\hang{\bf Corollary 6:}  {\it If the metric on $\pi^{-1}(b)$
is flat and if the connection on $E$ is flat in fiber directions
in $\pi^{-1}(b)$, then $R_1^Q =\str^Q(\Omega^u)$ vanishes at $b$.}
\bigskip

This follows from (6.12), since $\alpha_0(x)\equiv 0$ under the
hypotheses.  The significance of the Corollary is that
the curvature of the Bismut-Freed connection on the determinant line
bundle of a family of Dirac operators, given by
$$\Omega^{\Det, 0} = \biggl[\int_M \hat
A(\Omega^u)\ch(\Omega^E)\biggr]_{[2]},$$ 
does not have this vanishing property.  Thus Theorem 3 splits the
Bismut-Freed
 curvature into two terms: $R_1^Q = \str^Q(\Omega^u)$, 
the analogue of the finite
dimensional curvature, is localized on
the fiber in the sense of the Corollary; the other obstruction term is
a non-local Wodzicki residue.

   \bp $\bullet$ {\bf  Remarks:} 
(1) We can define a Chern character form as $\sum_k 
 \str^Q(\Omega^k)/k!$ for weighted bundles, and hence Chern forms
 via Newton's formulas.  These forms will not be closed in general,
 and their significance is unclear. 

(2) Theorems 3 and 5
compute the infinite dimensional obstruction to the finite dimensional
equality of 
  the curvature on the determinant bundle 
with (minus) the  first Chern form on the original vector bundle. 
The different looking obstructions in these theorems 
are
   related by the fact that renormalized limits of expressions of the type
$\langle  A_0, A_1, \cdots,
 A_k\rangle_{k, \e, Q_0}$ can be expressed in terms of Wodzicki residues.
More precisely, in  Appendix C we show that the
%terms of the type $\Lim^{\mu}\left(   \e \langle I, [\nabla^{\cal E}, L],
%[\nabla^{\cal E}, L]\rangle_{\e, 2, Q_0}\right)$ 
coefficients
of divergent terms
 in the asymptotics
 of  $    \langle I, [\nabla^{\cal E}, L],
[\nabla^{\cal E}, L]\rangle_{\e, 2, Q_0}  $ as $\e \to 0$ 
are combinations of
Wodzicki residues. 

(3) In fact, Proposition 4
is a more refined  result than Theorems 3 and 5.
Indeed,   zeta function regularization only
detects logarithmic divergences, while   heat kernel
regularization  
keeps track of
all  divergences in
fractional powers of $\e.$

\bigskip

\noindent {\twelvepoint {\bf 7.
 The Bismut-Freed connection and the
 curvature of the  determinant bundle over the manifold of almost complex
structures  }}

\tenpoint

%\noindent \centerline{ \bf  Appendix B:
% The Bismut-Freed connection and the
% curvature of the  determinant bundle over the manifold of almost complex
%structures }

   \np In this section, we apply the theory of \S6 to
study the Bismut-Freed connection
   on the fibration associated to the string theory example of
   diffeomorphisms acting on the space of almost complex structures 
${\cal A}(\Lambda)$ on
   a surface.  We show that this connection agrees with a
   classical connection in Teichm\"uller theory, and we compute
   the renormalized first Chern form for the infinite dimensional bundle.

Let 
 $\Lambda$ be a smooth
 closed Riemannian surface of genus greater than one,  and fix a
 Sobolev index   $s>1$.
As in \S1,     Example iii), we set
$${\cal E}^+ = {\cal A}(\Lambda)\times H^{s+1}(T\Lambda), $$
 $$  {\cal E}^- = T{\cal A}^s\mid_{{\cal A} (\Lambda)}=
\bigcup_{J\in {\cal A}(\Lambda)}
 \{H\in H^s(T^1_1\Lambda), J H=-HJ\}. $$

  \np $\bullet$ {\bf
 Almost complex structures on the bundles ${\cal E}^\pm$:}  
Each of the real bundles ${\cal E}^\pm$  over
${\cal A}(\Lambda)$  has an almost complex structure. On ${\cal E}^+$,
the almost complex structure  is  defined on the fiber
 above $J\in {\cal A}(\Lambda)$
  by $J$ itself:
 $${\cal J}^+_J(u)\equiv  J u,\ u\in
H^{s+1}(T\Lambda).$$
Similarly, the action
$${\cal J}^-_J(H)= H\cdot J,\   H\in T_J {\cal
A}^s(\Lambda)$$
is an almost complex structure on ${\cal E}^-$.

\np Let ${\cal M}^s_{-1}(\Lambda)$ be the space of $H^s$ Riemmanian metrics on $\Lambda$ with curvature $-1$, and set
 $$  \Phi: {\cal A}^s(\Lambda)\to   {\cal
M}_{-1}^s(\Lambda),\ \Phi(
J) = g_J, \eqno(7. 1)$$ 
where $g_J$ is the 
 unique Riemannian metric $g_J $  on $\Lambda$ with curavutre $-1$ in
 the conformal class defined by $J$. 
 $\Phi$  is a diffeomorphism  between
the Hilbert manifolds
 ${\cal A}^s(\Lambda)$ and ${\cal M}_{-1}^s(\Lambda)$, and the
derivative of  $\Phi$  at
$J$ in the direction $N$ is given by
$\left( D_J\Phi(N)\right)_{ab} ={ (g_J )}_{ac }(N  J)^c_b $ [T].
For $\alpha_J$ as in (2.1), the operator 
$$P_g = P_{g_J}\equiv   D_J \Phi  \circ  \alpha_J \circ
 D_{g_J }\Phi^{-1} \eqno(7. 2)$$
plays  a fundamental role in the
Faddeev-Popov procedure for string theories (see
[AJPS]).

\bigskip
\noindent\hang{\bf Lemma 7:}
{\it     The bundle map $\alpha :{\cal E}^+ \to {\cal E}^-
  $ defined in (2.1) is compatible with the
 almost complex structures ${\cal J}^\pm$ in the sense that}
$$\alpha_J( \Ker( {\cal J}_J^+-i))= \Ker ({\cal J}_J^--i),\quad
\alpha_J( \Ker( {\cal J}_J^-+i))= \Ker ({\cal J}_J^-+i),\  \forall J\in
{\cal A}(\Lambda).$$
{\it  Moreover,  $\alpha_J$ is a first order
elliptic operator.}
\bigskip

\noindent {\bf Proof: } 
We first show that   $\alpha_J$ is first order elliptic.
In isothermal coordinates for  $g$, the complexified operator
$P_g^{\C}$ is
$$  P_g^{\C}(u^{\bar z} {\partial \over \partial \bar z}+ u^z  {\partial
\over \partial  z}) =
   \partial_{ z }  u^{\bar z} {\partial \over \partial \bar z}\otimes dz  +
 \partial_{\bar z }  u^{  z} {\partial \over \partial  z}\otimes d \bar z
\eqno (7.3)$$
[AJPS], so $P_g^{\C}= \bar \partial \oplus \partial$ where $\bar \partial $
is the Cauchy-Riemann operator.
  $P_{g_J}$ is therefore first order elliptic, and hence so is
$\alpha_J$, since
 its principal symbol differs from $P_{g_J}$'s by the
isomorphisms in
 (7.2).

 \sp It is easy to check that
$\Ker({\cal J}_J^+-i)=
\{u^z {\partial \over\partial z}\}$  and  $\Ker( { \cal J}_J^++i)=
\{u^{\bar z} {\partial \over\partial\bar z}\},
$  and that
$$\Ker({\cal J}_J^--i)=
\{H_{\bar z}^z {\partial \over\partial  z}d\bar z\},\ \Ker({\cal
J}_J^-+i)=
\{H^{\bar z}_z {\partial \over\partial \bar z}dz\}.$$ 
Indeed, since $J^1_1=
J^2_2=0$ and $J^1_2=-1$ in
isothermal coordinates,
 we have
$$\eqalign{ ({\cal J}_J^-H)_z^{\bar z} &= (HJ)_z^{\bar z}= {1\over 2}
((HJ)_2^1- i (HJ)_1^1)
= {1\over 2} (H_1^1 J_2^1- i H_2^1 J_1^2)\cr
&= {1\over 2} (-H^1_1-i H^1_2) 
= -{i\over 2} (H^1_2-i H^1_1)
=-iH_z^{\bar z}.\cr}$$ 
Similarly,
$ ({\cal J}_J^+H)^z_{\bar z}= i H^z_{\bar z}.$
%Thus  we have that $Ker({\cal J}_J^-+i)= \{ H_z^{\bar z} {\partial \over
%\partial \bar z}
%\otimes d  z\}$ and
%  $Ker({\cal J}_J^--i)= \{ H^z_{\bar z} {\partial \over \partial   z}
%\otimes d \bar z\}$.
 The lemma then follows from (7.3), since
$$u\in \Ker ({\cal J}_J^+-i)\Rightarrow u^{\bar z}=0\Rightarrow
\left(P_{g_J} (u)\right)_z^{\bar z}=0\Rightarrow
P_{g_J} (u)\in \Ker({\cal J}_J^--i).$$
${}$\hfill$\bullet$

  \np $\bullet$ {\bf Hermitian metrics on ${\cal E}^\pm$:}
 ${\cal E}^\pm$  have the $L^2$
Riemannian metrics  $\gamma^\pm$ given by
(2.2${}^\pm$), which are compatible with the almost
complex structures
 ${\cal J}^\pm$. Indeed, for tangent vector fields
$u,v$
on $\Lambda$, we have
$$ \langle {\cal J}^+ u, {\cal J}^+ v\rangle_J^+=  \langle J u, J
v\rangle_J^+
= \langle u,   v\rangle_J^+,$$
since 
 $g_J$ is  compatible with  $J$.
Similarly, for $(1,1)$ tensors $H,K$ on $\Lambda$, we have
  $$ \langle {\cal J}^- H, {\cal J}^- K\rangle_J^-=  \langle   H
J, K J  \rangle_J^+
= \int_\Lambda d \mu_J(x) \tr (HJJ^*K^*)
= \int_\Lambda d \mu_J(x) \tr (H K^*),$$
since $J^*=-J$ and $J^2=-1$.

\sp Using the family of elliptic operators  $\{Q_J\equiv Q_J^+\oplus
Q_J^-\equiv 
 \alpha_J^* \alpha_J\oplus \alpha_J^* \alpha_J, 
\ J\in {\cal A}(\Lambda)\},$
we have $H^s$ metrics $\gamma^{s,\pm}$ 
 defined on the fiber ${\cal E}_J^\pm$
above
$J$ by
$$\langle u, v\rangle_J^{s,\pm}\equiv  \langle (Q_J^\pm +1)^{s} u, 
v\rangle_J^\pm=
\langle (Q_J^\pm+1)^{s\over 2} u,(Q_J^\pm+1)^{s\over 2} v\rangle_J^\pm.
\eqno(7.4)$$

 \np  $\bullet$ {\bf Connections on ${\cal E}^\pm$:}
We now define $L^2$ and $H^s$ connections on ${\cal E}^\pm.$
 As ${\cal E}^+$ is trivial,   let
$\nabla^+\equiv  d + \theta^+$ where
$$\theta^+(N) \equiv   {1\over 2} NJ,  \eqno(7.5^+)$$
for 
 $N\in T_J{\cal A} (\Lambda)$, $J\in {\cal A}(\Lambda)$. Here
 $NJ$ acts on  $u\in
H^{s+1}(T\Lambda)$ by
 $NJ u(x)\equiv  N(x) \cdot
J(x) \left(u(x)\right)$, with $\cdot$ denoting  matrix
multiplication. Since $N,J\in
\Ci(T^1_1\Lambda)$, multiplication by $NJ$ preserves
$H^{s+1}(T\Lambda)$, so  $\theta^+$ is
 a $\Hom(H^{s+1}(T\Lambda),H^{s+1}(T\Lambda))$-valued one-form 
on ${\cal A}$.

  \sp   The local charts
on the manifolds ${\cal A}^s(\Lambda)$  and ${\cal A} (\Lambda)$,
given pointwise by
the matrix exponential map as in \S1, 
 induce  a local trivialization of
${\cal E}^-$  over the
  base space ${\cal A} (\Lambda)$   with  fibers
 $T_J {\cal A}^s(\Lambda) $, $J\in {\cal A}(\Lambda)$. In a local
chart at
$J\in {\cal A} (\Lambda)$, we set
$\nabla^-\equiv  d + \theta^-$, with
$$\theta^-(N )\equiv   -{1\over 2}J\{ N, \cdot\},  \eqno(7.5^-)$$
for $N\in T_J{\cal A}(\Lambda).$
Here 
 $\{ M, N\}\equiv  MN+ NM$. 
 $\theta^-$ is a $\Hom(H^s(T^1_1\Lambda),H^s(T^1_1\Lambda))$-valued
 one-form on ${\cal
A}(\Lambda)$, since 
we again matrix multiply 
elements in $H^s(T_1^1\Lambda)$ by  elements in $\Ci(T_1^1\Lambda)$.  This
connection corresponds to
 the ``algebraic connection'' defined in [T, (5.6)].

\bigskip
\noindent\hang {\bf Lemma 8: } 
  1)\ {\it   $\nabla^\pm$ are  compatible
with the $L^2$-metrics $\gamma^\pm$ and with the almost
complex structures
 ${\cal J}^\pm$
in horizontal directions. In other words,  the L$^2$ 
superconnection $\nabla$ is
K\"ahler
 in horizontal directions.}

 2) {\it  $\nabla^{s, \pm}\equiv \left( Q^\pm+I\right)^{-{s\over
2}}\nabla^{\pm}
\left( Q^\pm+I\right)^{ s\over 2}$
 are  compatible
with the $H^s$-metrics $\gamma^{s,\pm}$ and
with the almost
complex structures ${\cal J}^\pm$ and ${\cal J}^-$
in horizontal directions. In particular, the connection
 $$\nabla^s\equiv (\nabla^{s, +})^*
\otimes 1+ 1\otimes \nabla^{s, -}$$ 
is
K\"ahler
 in horizontal directions.}
\bigskip

\noindent {\bf Remark:} It is shown in [T, Thm.~5.2.2] that in
horizontal directions,
  $\nabla^-$ equals
the $L^2$-Levi-Civita connection on the manifold of almost complex
structures.
 \bigskip

\noindent {\bf Proof: }  1)
The
compatibility of $\nabla^-$ with ${\cal J}^-, \gamma^-$ is shown in
[T, Thms.~5.2.1, 5.2.2].  We adapt this proof to $\nabla^+$ and refer
the reader to [T] for details.

\sp To prove the compatibility of
$\nabla^+$ with $\gamma^+$,  first note that
 the derivative of the map $g\mapsto \mu_g$ sending a Riemannian metric $g$
on $\Lambda$ to the
 corresponding volume form $\mu_g$  vanishes in the 
direction of a 
traceless covariant two tensor. Indeed, we have
$D_h (\mu_g)= {1\over 2} \tr_g (h) \mu_g= 0.$
 For any horizontal vector field $N$ at  $J$,
we
  set $n\equiv D_J\Phi(N) $, for $\Phi$ in (7.1),  $n$  is  a
traceless covariant two tensor
 [T, Thm.~2.5.6]. For $u,v\in\Gamma({\cal E}^+),$ we
have
$$\eqalign{ N \langle u, v\rangle_J^+
&=n\int_\Lambda  d\mu_{g_J}   {g_J}_{ab}u^a v^b
 \cr
&=\int_\Lambda    d\mu_{g_J}   n_{ab}   u^a v^b+ \int_\Lambda    d\mu_{g_J}
{g_J}_{ab}D_n u^a v^b
+ \int_\Lambda    d\mu_{g_J}    {g_J}_{ab}  u^a D_n v^b\cr
&=\int_\Lambda    d\mu_{g_J}   {g_J}_{ac} N^c_d J^d_b   u^a v^b+
\langle D_N u, v\rangle_J^+ + \langle u, D_N v\rangle_J^+ \cr
&=\langle        N J    u, v \rangle_J^+ +
\langle D_N u, v\rangle_J^+ + \langle u, D_N v\rangle_J^+ \cr
&={1\over 2}\langle        N J    u, v \rangle_J^+ -{1\over 2}\langle  u,
JNv \rangle_J^+ 
+
\langle D_N u, v\rangle_J^+ + \langle u, D_N v\rangle_J^+ \cr
&= \langle \nabla_N^+ u, v\rangle_J^+ + \langle u, \nabla_N^+  v\rangle_J^+,
\cr}\eqno (7.6)$$
where we have used  $JN=-NJ$ and $N^*=N$.

\sp For the compatibility with the almost complex structure
${\cal J}^+$, we have
$$\eqalign{ [\nabla^+_N, {\cal J}^+]u&= D_N(J u)+ {1\over 2} N J^2 u- J D_N
u- {1\over 2} J N J \cr
&=  N u -{1\over 2} N u+ {1\over 2} J^2 N u \cr
&= 0.\cr}$$

2) This is a straightforward
consequence of 1), once we check
  the  compatibility  of  $Q_J^\pm$
with
the almost complex structures ${\cal J}_J^\pm$, 
 which follows from    Lemma 7. \hfill $\bullet$

 \np
$\bullet$ {\bf A half-weighted vector bundle:} 
The   bundle
$$\Hom\left({\cal E}^{+,1,0} , {\cal E}^{-,1,0} \right)\equiv  \left({\cal
E}^{+,1,0}\right)^* \otimes
 {\cal E}^{-,1,0}  $$ 
now has a
connection $\nabla^s\equiv  {(\nabla^{s,+})}^*\otimes 1+ 1\otimes \nabla^{s,-}$
which is horizontally K\"ahler.
In a local chart, we have
$\nabla^s= d + \theta^s= d +\theta^{s,-}-\theta^{s,+}$, so
we can equivalently view $\nabla^s$ 
as a
superconnection on the superbundle
$${\cal E}^{1,0}\equiv  {\cal E}^{+,1,0} \oplus {\cal E}^{-,1,0}.\eqno(7. 7)$$
The family $J\to L_J^{1,0}\equiv  \left[\matrix { 0&L_J^-\equiv \partial_J\cr
L_J^+\equiv \bar \partial_J &0 \cr}\right],$
  where $\bar \partial_J$ is the Cauchy-Riemann operator for 
$(\Lambda, J)$,
defines a section of the bundle
$Ell\left({\cal E}^{  1,0}     \right)$. 
By Lemma 8,
 $  L_J^{1,0}$ is a  self-adjoint elliptic operator for the   hermitian product
built from the almost complex structure
${\cal J}\equiv  {\cal J}^+ \oplus {\cal J}^-$ and the scalar product
 $\langle \cdot , \cdot \rangle_J^+\oplus   \langle \cdot , \cdot
\rangle_J^- $.
 (cf.~[AJPS], [T] for a string theory perspective).  Thus
 $\left( {\cal E}^{1,0}, L^{1,0}\right)$ is a half-weighted vector
bundle.
$Q^{\pm,1,0}\equiv 
L^\mp L^\pm$ are positive
self-adjoint sections of
$Ell({\cal E}^{\pm, 1,0})$.  Here  $L^\mp$ is
either the  $L^2$ or
the
$H^s$ adjoint of $L^\pm$
with respect to the inner products (7.5${}^\pm$).  This data determines
a weighted superbundle
 $  \left({\cal E}^{1,0} ,Q^{1,0}\equiv 
 Q^{+,1,0}\oplus Q^{-,1,0}\right)$.

     \bp $\bullet$ {\bf The first Chern forms  of ${\cal E}^{\pm,1,0}$:}

 \noindent\hang {\bf Lemma 9:} {\it  Let $({\cal E}, Q)$ be a weighted vector
bundle with an almost
 complex structure ${\cal J}$ compatible with $Q$ (i.e.~$Q{\cal J} =
{\cal J} Q$),
let $({\cal E}^{1,0}, Q^{1,0})$ denote its $(1,0)$ part, and let
$A\in\Gamma(
PDO({\cal E}))$ satisfy $A{\cal J}= {\cal J} A$.  Then}
$$\tr^{Q^{1,0}} (A^{1,0})= \tr^Q(A)+ i\  \tr^Q({\cal J} A).$$ 
\bigskip

\noindent {\bf Proof: } Let $A^{1,0}$ be the $(1,0)$ part of
$A^{\C}$, the fiberwise complexification of $A$ with respect to
${\cal J}$.       It is standard that
$\tr(A^{1,0})= \tr(A)+ i\ \tr(JA)$.  Then
$$\eqalign{ \tr^{Q^{1,0}} (A^{1,0})&= \Lim_{z\to
0}\tr(A^{1,0} {(Q^{1,0})}^{-z})
=\Lim_{z\to 0}\tr({(A   Q^{-z})}^{1,0}) \cr
 &=\Lim_{z\to 0}\left( \tr(  A   Q^{-z}  ) + i \Lim_{z\to 0} \
\tr(  {\cal J} A
Q^{-z}  ) \right)
=   \tr^Q (  A    ) + i \tr^Q(  {\cal J} A  ). \cr}$$
${}$\hfill$\bullet$

\np  We now
compute the curvature of $\nabla^{  s,\pm}$ on ${\cal E}^\pm$.

\bigskip
\noindent\hang{\bf Lemma 10:}
{\it For $s> 1$, the curvatures $\Omega^{s, \pm}$,
 of the connections
$\nabla^{s, \pm}$  are zero order PDOs given by}
 $$\eqalign{ \Omega^{s, +}   (M, N)H
&= -{1\over 4} \left( Q^++I\right)^{-{s\over 2}}
 [M, N]_{\op} H
\left( Q^++I\right)^{ s\over 2},\cr 
  \Omega^{s, - }   (M, N)H &= \left( Q^-+I\right)^{-{s\over 2} }
\left( -{1\over 2} [M, N]_{\op}H   +{1\over 2}(-MHN+NHM)\right.\cr
&\qquad \left.
-{1\over 4}  \left[
\{M, H\}, \{N, H\}\right]\right)
\left( Q^- +I\right)^{ s\over 2},\cr}$$
{\it for $M,N,H\in T_J{\cal A}(\Lambda).$}
\bigskip

Here $[M, N]_{\op}$ denotes the multiplication operator in the fiber
 over $J$ associated to 
 the bracket (pointwise over $\Lambda$)
of the matrices $M,N.$  In contrast,  $[M,N]$ denotes
 the bracket of vector fields on ${\cal A}$ which are
 given by local extensions of the tangent vectors $M,N$ at $J$.  At a
 fixed $J$, we may extend $M,N$ so that $[M,N]=0.$
\bigskip

\noindent {\bf Proof:}   We prove the first equality only,
 since the  second is
similar.  Using $M(J) = M, JM= -MJ, J^2 = -{\rm Id}$ and similar
formulas for $N$,
we have
 $$\eqalign{ \Omega^{s, +}   (M, N)&=  \left( \nabla^{s, +}\right)^2  (M, N)
=[\nabla_M^{s,+},\nabla_N^{s,+}] - \nabla_{[M,N]}^{s,+}\cr
&= \left( Q^++I\right)^{-{s\over 2}}
 \left( [\nabla^{  +}_M, \nabla^{ +}_N] -
\nabla^{  +}_{[M, N]}\right)
\left( Q^++I\right)^{ s\over 2} \cr
&=\left( Q^++I\right)^{-{s\over 2}} \left(d \theta^+(M, N) + \theta^+ \wedge
\theta^+ (M, N) \right)
\left( Q^++I\right)^{ s\over 2} \cr
&=\left( Q^++I\right)^{-{s\over 2}} \left( M( \theta^+(N))- N( \theta^+(M))
-\theta^+([M,N]) + \theta^+ \wedge
\theta^+ (M, N) \right)
\left( Q^++I\right)^{ s\over 2} \cr
&= \left( Q^++I\right)^{-{s\over 2}} \left(-{1\over 2} [M, N]_{op} + {1\over
4}[MJ, NJ]\right)
\left( Q^++I\right)^{ s\over 2}\cr
&= \left( Q^++I\right)^{-{s\over 2}} \left( -{1\over 4} [M, N]_{op} \right)
\left( Q^++I\right)^{ s\over 2}.\cr  }$$
${}$\hfill$\bullet$

 \bigskip
\noindent\hang {\bf Proposition 11:} {\it   
The weighted first Chern form $R_1^{Q
}$ on the weighted vector bundle
$\left({\cal E}^{1,0}, Q^{1,0}\right) $  with the connection $\nabla^s$
is independent of the parameter $\mu$ used in the renormalization
procedure and 
independent of $s>1$. For  $M, N\in T_J{\cal A}(\Lambda)$, we have}
$$\eqalign{R_1^{Q }(M,N) &= 
  i\ \tr^{Q^- }  \left( {1\over 2}J    [M, N]_{\op} 
-{1\over 2}(-M(\cdot)N + N(\cdot)M) 
+{1\over 4}  J  \left[
\{M, \cdot\}, \{N, \cdot\}\right]  \right)\cr
&\qquad +{i\over 4} \tr^{Q^+ }  ( J    [M,
N]_{\op}   ).\cr}$$

\bigskip

 The traces are taken with respect to the $L^2$ inner products.  Note
that, in agreement with Corollary 6, the curvature on the associated
determinant bundle is the $Q$-weighted trace of a multiplication
operator.
\bigskip

\noindent {\bf Proof:}    By Lemma 7,  $Q^+$ commutes with the
almost complex structure, so the     $(1,0)$ part
 $\Omega^{s, +,1,0  }   (M, N)
 $ of   $ \Omega^{s, +  }   (M, N)$ satisfies
 $\Omega^{s, +,1,0  }   (M, N)= \left(\Omega^{s, +}    (M,
N)\right)^{1,0}.$ 
Applying Lemmas 9 and 10, we find
 $$ \eqalign{  \tr_s^{{Q^+}^{1,0}, \mu}  (\Omega^{s, +,1,0 }   (M, N)  ) &=
  \tr_s^{Q^+, \mu }  (\Omega^{s, +  }   (M, N)  ) + i\ 
\tr_s^{Q^+ , \mu}  ({\cal
J}^+_J\Omega^{s, +  }   (M, N)  )\cr
&=  -{i\over 4} \tr^{Q^+, \mu }  ( J    [M, N]_{\op}   ).\cr}$$
Note that $\tr_s^{{Q^+}, \mu}  (\Omega^{s, + }   (M, N)  )
=0$, since the curvature form is a skew-symmetric endomorphism for
fixed $M,N$.  
The $H^s$ trace $\tr_s$ defined via  the inner
 product $\langle \cdot, \cdot\rangle_J^{s,+}$ in the first line equals
the
$L^2$ trace in the second line, because the powers $\left( Q^+\right)^
{\pm{s\over
2}}$
cancel in the computation of $\tr_s.$  
Similarly, we obtain
 $$ \eqalign{  \tr_s^{{Q^-}^{1,0}, \mu}  (\Omega^{s, -,1,0 }   (M, N)  ) &=
  \tr_s^{Q^-, \mu }  (\Omega^{s, -  }   (M, N)  ) + i\ \tr_s^{Q^- , \mu}  ({\cal
J}^-_J\Omega^{s, -  }   (M, N)  )\cr
&= i\   \tr^{Q^-, \mu}  \left( {1\over 2}J    [M, N]_{\op}   
-{1\over 2}(-M(\cdot)N + N(\cdot)M) 
+{1\over 4} ( J  \left[
\{M, \cdot\}, \{N, \cdot\}\right]  \right).\cr}$$
For a differential operator $A$, there is no  logarithmic divergence
in the asymptotics of $\tr(A e^{-\e Q})$ as $\e \to 0$.  Since $\mu$ keeps track of the logarithmic divergence,
the renormalized traces above are independent of 
$\mu$.

\sp   The result now  follows from Definition (4.8).
\hfill$\bullet$ 
\bigskip

\sp\noindent {\bf Remark:} 
 A matrix $H\in T_J{\cal A}^s(\Lambda)$ satisfies  $HJ=-JH$ so
 two such matrices
$H, K$ satisfy $HKJ=JHK$. Writing $J=\left[\matrix{ 0 &-1\cr
1& 0\cr}\right]$ in isothermal coordinates, we see that $HK$ is of the form
$\left[\matrix{ \alpha &\beta\cr
\beta & -\alpha\cr}\right],$ as is any even product of matrices
in $T_J{\cal A}^s(\Lambda)$. Hence  $J    [M, N]_{\op}  $ is of the form
 $\left[\matrix{ \gamma &\delta\cr
-\delta&  \gamma\cr}\right]$. In contrast to an incorrect claim in
 [PR], 
$\tr^{Q^+,\mu }  ( J    [M, N]_{\op}   )$ need not vanish.

\bp\noindent $\bullet$ {\bf $\nabla^+$ and the Bismut-Freed
 connection:}  We now show that the connection $\nabla^+$ of (7.5${}^+$) coincides with the
Bismut-Freed connection associated to the string theory fibration
$\Lambda\to ({\cal A}\times \Lambda)/{\cal D} \to {\cal A}
/{\cal D}$, where ${\cal D} = {\rm Diff}_0^{s+1}$ is the (Sobolev
$s+1$) isotopy group 
 of $\Lambda.$  It is equivalent to work with the
trivial fibration $\Lambda \to {\cal A}\times \Lambda\to {\cal A}$,
and to consider only directions perpendicular to the action of ${\cal
D}$ on ${\cal A}\times \Lambda$ with respect to the natural metric
$$\Vert (h,v)\Vert_x^2 = \Vert h\Vert_{g_J}^2 + \Vert
v\Vert^2_{g_J},\eqno (7.8)$$ 
where $h\in T{\cal A}, v\in T\Lambda$, and the projection
$\pi:{\cal A}\times \Lambda\to{\cal A}$ has $\pi(x) =J.$ Note
that the role of $\nabla^-$ of (7.5${}^-$) is implicit, since it is the Levi-Civita
connection for the metric on $T{\cal A}$.

\sp Let $\nabla^{\LC}$ be the Levi-Civita connection on ${\cal A}\times
\Lambda$ for the metric (7.8).  By (6.10), 
for $\psi\in\Gamma(T\Lambda)$, we must show that
$$\nabla_h^+\psi = P^{T\Lambda}\nabla_h^{\LC}\psi,\eqno(7.9)$$
where $P^{T\Lambda}$ is the orthogonal projection of $T({\cal A}\times
\Lambda)$ to $T\Lambda.$  By the six term formula for the Levi-Civita
connection, we have
$$\eqalign{ 2\langle P^{T\Lambda}\nabla_h^{\LC}\psi,v\rangle &= 
2\langle \nabla_h^{\LC}\psi,v\rangle\cr
&= h\langle \psi,v\rangle +\psi\langle h,v\rangle -v\langle
h,\psi\rangle\cr
&\qquad +\langle [h,\psi],v\rangle + \langle [v,h],\psi\rangle -
\langle [\psi,v],h\rangle.\cr}\eqno(7.10)$$
On the right hand side of (7.10), we may extend $h,v$ arbitrarily near
$x$, so we choose $h$ to be horizontal and $v$ to be vertical near
$x$.  Then 
$$\langle h,v\rangle = \langle h,\psi\rangle =0\eqno (7.11)$$
 in (7.10).

\sp Let $\phi_t^v$ be the vertical flow of $v$.  Then
$$\eqalign{\langle [v,h],\psi\rangle &= \langle {d\over dt}\biggl|_{_{t=0}}
\phi^v_{-t,*} h,\psi\rangle_{g_J} = {d\over dt}\biggl|_{_{t=0}}\langle
\phi^v_{-t,*} h,\psi\rangle_{g_J} \cr
&= \langle h,[-v,\psi]\rangle = \langle [\psi,v], h\rangle,\cr}\eqno(7.12)$$
(Since ${\cal D}$ preserves the volume measure $d\mu_{g_J}$, we may
move ${d/dt}$ past the inner product in the first line.)  Combining
(7.10)-(7.12) gives
$$\eqalign{2\langle P^{T\Lambda}\nabla_h^{\LC}\psi,v\rangle &= 
h\langle \psi,v\rangle + \langle [h,\psi], v\rangle\cr
&= \langle \nabla_h^+\psi,v\rangle +\langle \psi,\nabla^+v\rangle 
+ \langle [h,\psi], v\rangle\cr
&=\langle h(\psi),v\rangle + \langle \psi,h(v)\rangle + \langle
[h,\psi], v\rangle,\cr}$$
where we have used the third line of (7.6) in the last line.
Moreover, $[h,\psi] = h(\psi) - \psi(h) = h(\psi)$, since $h\in T{\cal
A}$ may be lifted to be constant in vertical directions.  So we
finally obtain
$$\eqalign{ 2\langle P^{T\Lambda}\nabla_h^{\LC}\psi,v\rangle &=
2\langle h(\psi),v\rangle + \langle hJ\psi,v\rangle + \langle \psi,
h(v)\rangle \cr
&= 2\langle h(\psi),v\rangle + \langle hJ\psi,v\rangle\cr
&= 2\langle\nabla_h^+\psi,v\rangle,\cr}$$
since the extension of $v$ may be taken to be constant in vertical
directions, and so $h(v)=0.$

   \vfill \eject

\noindent {\twelvepoint {\bf Appendix A: Superconnection formalism  }}

\tenpoint
%\noindent \centerline{ \bf  Appendix A: Superconnection formalism
%}\np

\np
This appendix summarizes the superconnection formalism used in \S5.
Useful references are [BGV], [Q2].

  \np $\bullet$ {\bf Super vector bundle valued forms:} A super vector bundle
  ${\cal E}$ over a manifold $B$ is a $\Z_2$ graded vector bundle
  ${\cal E}=  {\cal E}^+ \oplus {\cal E}^-$ over $B$. Let $\Omega(B ,
{\cal E})$ be
  the space of ${\cal E}$-valued differential forms on $B$.  Since
  $\Omega(B , {\cal E})\simeq \Omega(B)\otimes_{C^\infty(B)}  \Gamma({\cal
  E}) $, where
  $\Omega(B)$ is the exterior algebra of forms on $B$, 
 an element of $\Omega(B, {\cal E})$ can be written as
$\omega \otimes \sigma$ with $\omega\in \Omega(B),\
  \sigma\in\Gamma({\cal E}).$
  The $\Z $ grading on $\Omega(B)$ induces a $\Z_2$ grading on
  $\Omega(B)= \Omega^+(B) \oplus \Omega^-(B)$ into forms of even and odd degree,
which, with the
$\Z_2$ grading
  on ${\cal E}$, yields a
  $\Z_2$ grading on $\Omega(B, {\cal E}) = \Omega^+(B, {\cal E})\oplus
  \Omega^-(B, {\cal E})$, where
  $$\eqalign{\Omega^+(B, {\cal E}) &\equiv  \Omega^{+} (B, {\cal
E}^+) \oplus
   \Omega^{-} (B, {\cal E}^-),\cr
\Omega^-(B, {\cal E}) &\equiv  
\Omega^{+} (B, {\cal E}^-) \oplus
 \Omega^{-} (B, {\cal E}^+).\cr}$$

\np $\bullet$ {\bf From a connection to a one-parameter family of
 superconnections:} 
A superconnection is an odd
 first order differential operator:
$\nabla: \Omega^\pm(B, {\cal E})\to \Omega^\mp (B, {\cal E})$ which satisfies the
Leibniz  rule in the $\Z_2$ graded sense:
$$\nabla (\omega\otimes \sigma)= d  \omega \otimes \sigma + (-1)^{\vert
\omega\vert}
\omega \otimes \nabla \sigma.$$
A connection $\nabla$ on ${\cal E}$ which preserves the $\Z_2$ grading 
defines a map
$$\nabla: \Gamma({\cal E}^\pm) \left(\subset \Omega^\pm(B, {\cal E})\right)
\to \Gamma(T^*B \otimes {\cal E}^\pm)\subset \Omega^\mp (B, {\cal
E}),$$
which extends uniquely to a superconnection on ${\cal E}$.

\sp
The $\Z_2$ grading on ${\cal E}$ induces  a  $\Z_2$ grading on the
bundle
$\Hom(
{\cal E},{\cal E})= \Hom^+({\cal E})\oplus \Hom^-({\cal E})$, where the even
bundle maps, the sections of
$\Hom^+({\cal E})$, preserve the $\Z_2$ grading on
${\cal E}$, and the odd bundle maps, the sections of
$\Hom^-({\cal E})$, take ${\cal E}^\pm$
to 
${\cal E}^\mp$. A section $L$ of $\Hom^-({\cal
E})$ induces an odd map
$$L: \Omega^\pm(B, {\cal E})\to \Omega^\mp (B, {\cal E}),
\  \omega \otimes \sigma \mapsto (-1)^{\vert \omega\vert }\omega \otimes
L\sigma.$$
$\nabla$ and $L$ induce
a one-parameter family of
superconnections $\nabla_t^L\equiv  \nabla+ \sqrt t L, \ t>0$, on
${\cal E}$.

\np $\bullet$ {\bf From a superconnection on ${\cal E}$ to a
superconnection on $\Hom( {\cal E},{\cal E})$:} A superconnection $\nabla$ on
${\cal E}$ induces a connection on $\Hom({\cal E},{\cal E})$ defined by
$$[\nabla,A] \equiv  \nabla  A -(-1)^{\vert A\vert}A \nabla,
$$ 
where $\vert A\vert=0$ if $A$ is even and $\vert A\vert=1$ if
$A$ is odd.
If $A=L\in\Gamma(\Hom^-({\cal E}))$ and $\nabla $ is a superconnection induced by
a $\Z_2$ grading preserving  connection on ${\cal E}$,
then
$$\eqalign{ [\nabla, L](\omega \otimes \sigma)&\equiv  \nabla (L(\omega \otimes
\sigma)) +  L( \nabla (\omega \otimes \sigma))\cr
&=\nabla\left((-1)^{\vert \omega\vert }\omega \otimes L\sigma\right)+
L\left( d  \omega \otimes \sigma + (-1)^{\vert \omega\vert}
\omega \otimes \nabla \sigma\right) \cr
&= (-1)^{\vert \omega\vert }d\omega \otimes L\sigma+
(-1)^{\vert \omega\vert+1 }d\omega \otimes L\sigma\cr
&\qquad +(-1)^{2\vert \omega\vert}
\omega \otimes\nabla L\sigma+ (-1)^{2\vert \omega\vert+1 }\omega \otimes
L\nabla\sigma\cr
&=\omega \otimes [\nabla,L] \sigma.\cr}$$
In the last line, the bracket is an ordinary bracket.
 \vfill \eject

\noindent {\twelvepoint {\bf Appendix B: Trace forms and Wodzicki residues}}

\tenpoint

%  \centerline{ \bf Appendix C: Trace forms and Wodzicki residues}
\medskip 

In this appendix, we express  the divergences in the asymptotics of
the 
trace forms $\langle A_0, A_1, \ldots , A_k\rangle_{\e, k, Q}$ as $\e \to
0$ in terms of Wodzicki
 residues.  Such a relation is suggested by Theorems 1 and 3, which compute
the obstruction to the equality of the determinant curvature and the
renormalized first Chern form alternately as such a divergent term and
as a Wodzicki residue, respectively.
 Such trace forms have occurred in  quantum algebras studied
by  Jaffe, Lesniewski and Osterwlader [JLO] and in local index theory
in non-commutative geometry treated by Connes and Moscovici [CM].

 \smallskip \noindent {\bf Notation:}
 For $j\in\N$ and $A,Q
\in PDO(M, E)$ such that the $Q$ has scalar 
top order symbol,  $ [ A]_Q^j\in PDO(M, E)$ is the operator
 defined inductively by
$$ [ A]_Q^0\equiv  A, \quad [ A]_Q^{j+1}\equiv  [Q, [ A]_Q^j]. $$
We will often drop the subscript $Q$, and  use notation from
the body of the paper.
Notice that
the operator $[ A]_Q^j$ is of order $a+j(q-1)$ where $a=\ord(A),
\ q = \ord(Q)$,
and that $[ A]_{\e Q}^j= \e^j [ A]_Q^j$ for any $\e >0$, $j \in \N$.

\bigskip
\noindent\hang {\bf Lemma B.1:} [L, Lemma 4.2]
{  \it  If $p, \e, N>0$ satisfy ${N-a\over q}-p-\e>0$,
then
 $$e^{-t Q} A= \sum_{ j=0}^{N-1}{ (-t)^j \over j!}[A]_Q^j e^{-t Q} +
  R_N(A, Q, t),$$
where for any $c>0$ such that $Q+c$ is invertible, there exists $C>0$ such that
$\Vert R_N(A, Q, t) (Q+c)^p\Vert \leq C t^{{ N-a \over q}-p-\e}  $.}
\bigskip

\noindent\hang {\bf Lemma B.2:}
{\it Given  $A_0, A_1, \cdots, A_k \in PDO(M, E)$ and $j_k\leq N_k\in \N$,
 there exist
  $N_1 , N_2 ,\cdots,
  N_{k-1} \in \Z $ such that for $j_i\leq N_i$ and  $\alpha_i \in \{0,
1\}, i=1,
 \cdots,k$ with at least one $\alpha_i$ unequal to one,
 the operator $$\eqalign{& A_0 \left(R_{N_1}(A_1, Q,
\sigma_0)\right)^{1-\alpha_1}
 \left([A_1]_Q^{j_1}\right)^{ \alpha_1}
\cdot\ldots\cdot
 \bigl((R_{N_i}(A_i, Q,   \sigma_0+ \sum_{l=1}^i \alpha_l \sigma_l
)\bigr)^{1-\alpha_i}\left(
 [A_i]_Q^{j_i}\right)^{ \alpha_i}\cdot \ldots\cr
& \cdot
  \bigl(R_{N_k }(A_k, Q,   \sigma_0+ \sum_{i=1}^k \alpha_i \sigma_i
)\bigr)^{1-\alpha_k} \left([A_k]_Q^{j_k}
\right)^{\alpha_k}\cr}\eqno({\rm B}.1)$$
is trace-class with trace bounded by
$$C\cdot  \prod_{j=1}^k (\sigma_0+ \sum_{i=1}^j \alpha_i \sigma_i)^{ (1-\alpha_j)
\left( N_j- a_j \over q- p_j\right)  }$$
for some positive constant $C$.}
\bigskip

\noindent {\bf Proof:} We proceed by induction on $k$.
 For $k=1$, there is an integer
  $p_0$ such that
$A_0Q^{-p_0}$ is trace-class. By Lemma B.1,
we can choose $N_1$ such
 that  $ Q^{ p_0  } R_{N_1}(A_1, Q, \sigma_0)$ is bounded by 
$$C\sigma_0^{{N_1-a_1\over q}- p_0},\eqno({\rm B}.2)$$
 where $a_1= \ord(A_1)$ and $C$  is a positive constant. Then 
$A_0 R_{N_1}(A_1, Q, \sigma_0)$ is trace-class with trace bounded by an
expression similar to (B.2).

We now assume the  lemma through $k-1$ for the induction step. $C$ will
denote a constant which may change from line to line.
By Lemma B.1, there exists $N_k\in\Z$ such that
 $R_{N_k }(A_k, Q,   \sigma_0+ \sum_{i=1}^k \alpha_i \sigma_i)$ is bounded
by $C(\sigma_0+ \sum_{i=1}^k \alpha_i \sigma_i)^{N_k- a_k\over q}$.
For $\alpha_{k }=0$,  by induction  we can
choose
$N_1, \cdots, N_{k-1}$ such that
$$A_0 \left(R_{N_1}(A_1, Q,
  \sigma_0)\right)^{1-\alpha_1} \left([A_1]_Q^{j_1}\right)^{ \alpha_1}
\cdots
 R_{N_{k-1}}(A_{k-1}, Q,   \sigma_0+ \sum_{l=1}^{k-1} \alpha_l \sigma_l
) [A_{k-1}]_Q^{j_{k-1}}$$
 is trace-class with trace bounded by
$$C\cdot   \prod_{j=1}^{k-1}
(\sigma_0+ \sum_{i=1}^{j} \alpha_i \sigma_i)^{ (1-\alpha_j)
\left( N_j- a_j \over q- p_j\right)  }.$$
It follows 
from Lemma B.1 that
$R_{N_k}(A_k, A, \sigma_0+ \sum_{i=1}^k \alpha_i \sigma_i)$ is bounded in norm
by \hfill\break
  $C\cdot
( \sigma_0+ \sum_{i=1}^k \alpha_i \sigma_i)^{ N_k-a_k\over q}$. Hence (B.1)
is bounded
by 
$$C\cdot  \prod_{j=1}^{k-1} (\sigma_0+ \sum_{i=1}^j \alpha_i
 \sigma_i)^{ (1-\alpha_j)
\left( N_j- a_j \over q- p_j\right)  }\cdot
( \sigma_0+ \sum_{i=1}^k \alpha_i \sigma_i)^{ N_k-a_k\over q}. $$

\sp
Now assume  $\alpha_k=1$. We can choose
$p_{k-1},\ N_{k-1}$ large enough that
$Q^{-p_{k-1}}[A_k]^{j_k}_Q$ is bounded and
${N_{k-1}-a_{k-1}\over q}-p_{k-1}>0$. Then Lemma B.1 implies
$$\Vert R_{N_{k-1}}(A_{k_1}, Q,
 (\sigma_0+ \sum_{l=1}^{k-1 } \alpha_l \sigma_l ) Q^{p_{k-1}}\Vert\leq C
(\sigma_0+ \sum_{l=1}^{k-1 } \alpha_l \sigma_l)^{{N_{k-1}-a_{k-1}\over
q}-p_{k-1}}.$$
If $\alpha_{k-1} =0$, this estimate and the lemma for
$k-2$ produces the upper bound
$$C\cdot  \prod_{j=1}^{k-2} (\sigma_0+ \sum_{i=1}^j \alpha_i
 \sigma_i)^{ (1-\alpha_j)
\left( N_j- a_j \over q- p_j\right)  }\cdot
( \sigma_0+ \sum_{i=1}^{k-1} \alpha_i \sigma_i)^{{ N_{k-1}-a_{k-1}\over q}
-p_{k-1}}$$
for the trace.
If
 $\alpha_{k-1}=1$, there exist $p_{k-2},\ N_{k-2}$ such that
$Q^{-p_{k-2}} [A_{k-1}]_Q^{j_{k-1}} [A_k ]_Q^{j_k}$ is bounded and
 ${N_{k-2}-a_{k-2}\over q}-p_{k-2}>0$. 
Applying the above procedure gives the desired 
estimates.\hfill$\bullet$
${}$
\bigskip

 \noindent\hang{\bf
 Proposition B.3:}
{\it  Let $A_0, A_1, \cdots, A_k \in PDO(M, E)$. There exist
  $N_1 , N_2 ,\cdots,
  N_k\in\Z  $  such that for  $\e>0$}
$$\eqalign{ &\langle A_0, A_1, \cdots, A_k\rangle_{\sqrt \e, k,  Q}\cr
&\qquad=
\sum_{j_1=0}^{N_1-1} \cdots \sum_{j_n=0}^{N_k-1}{ \e^{j_1+ \cdots +j_n}
 \over j_1! \cdots j_k!} \int_{ 0}^1\cdots \int_{ 0}^1
(-1)^{j_1+\cdots j_k}(\sigma_0)^{j_1}\cdot  ( \sigma_0+ \sigma_1) ^{j_2}
 \cdots  (\sigma_{0 }+ \sigma_{ 1}+ \cdots\cr
 &\qquad\qquad+
\sigma_{k-1 })^{j_k}\tr(A_0
[A_1]_Q^{j_1}  \cdots [ A_k]_Q^{j_k}e^{- \e Q})d\sigma_0\cdots d\sigma_{k}
 +
{\rm o} (\e).\cr}$$

\bigskip

 \noindent
{\bf Proof: }Iterating Lemma B.1, 
we find
$$\eqalign{ &\langle A_0, A_1, \cdots, A_k\rangle_{  \e , k, Q }\cr
&\qquad =
\sum_{j_1=0}^{N_1-1} \cdots \sum_{j_n=0}^{N_k-1}(-1)^{j_1+\cdots j_k}
{ \e^{j_1+ \cdots +j_k}
 \over j_1! \cdots j_k!} \int_{ 0}^1\cdots \int_{ 0}^1
(\sigma_0)^{j_1}\cdot  ( \sigma_0+ \sigma_1) ^{j_2}
 \cdots  (\sigma_{0 }+ \sigma_{ 1}+ \cdots
+
\sigma_{k-1 })^{j_k}\cdot \cr 
&\qquad\qquad \cdot \tr(A_0
[A_1]_Q^{j_1}  \cdots [ A_k]_Q^{j_k}e^{- \e Q})d\sigma_0\cdots d\sigma_{k} +
R_{N_1, N_2, \cdots, N_k} (\e), \cr} $$
with
$R_{N_1, N_2, \cdots, N_k } (\e)$ a finite linear combination of
    terms of the type
$$\eqalign{ &\int_{ \sigma_0 + \cdots + \sigma_k=1, \sigma_i\geq 0}
\tr\biggl(A_0 (R_{N_1}(A_1, Q, \sigma_0))^{1-\alpha_1}
\left(\sum_{j_1=0}^{N_1} {( -\sigma_0 )^{j_1} \over j_1!}  [ A_1]_Q^{j_1}
\right)^{ \alpha_1} \cr
&\qquad \cdot (R_{N_2}(A_2, Q, \sigma_2))^{1-\alpha_2}
\left(\sum_{j_2=0}^{N_2} {(-(\sigma_0+ \alpha_1\sigma_1))^{j_2} \over j_2!}
{[ A_2]}_Q^{j_2}\right)^{ \alpha_2}\cdot
\cr
&\qquad\cdots \cr
&\qquad\cdot \left(R_{N_k}(A_n, Q, \sigma_k)\right)^{1-\alpha_k}
\left(\sum_{j_n=0}^{N_n} {(-(\sigma_0+ \alpha_1\sigma_1+ \cdots +
\alpha_k\sigma_{k }) )^{j_k}   [A_k]_Q^{j_k}  \over j_k!} \right)^{
\alpha_k} \cr
&\qquad \cdot e^{-\e (\sigma_0+ (1-\alpha_1)
\sigma_1+ \cdots +  (1-\alpha_k)\sigma_k) Q}\biggr)d \sigma_0 \cdots
d \sigma_{k },\cr}$$
with $\alpha_i$ equal to $0 $ or $ 1$,  and  $(\alpha_1, \cdots, \alpha_n)
\neq (1, \cdots, 1)$.
By Lemma B.2, $N_1, N_2, \cdots, N_k$ can be chosen so that
 the integrals 
in $R_{N_1, N_2, \cdots, N_k}(\e)$  converge and 
 $R_{N_1, N_2, \cdots, N_k}(\e)={\rm o}(\e)$. \hfill$\bullet$
 \bigskip 

\noindent
  $\bullet$ {  \bf  The asymptotics of regularized  trace  forms: }

 \smallskip  We now investigate the asymptotic behavior of the
 trace forms as $\e\to 0 $.
\bigskip

\noindent\hang       {\bf Theorem B.4:}
{\it Let  $A_0, \cdots, A_n \in PDO(M, E)$. Then
\item{(i)} $\langle A_0, A_1, \cdots, A_n\rangle_{  \e, n, Q}$ has the following
asymptotic expansion as $\e\to 0$:}
$$\eqalign{\langle A_0, A_1, \cdots, A_n\rangle_{  \e,n, Q}
&\sim \sum_{j=0}^\infty
 \alpha_{j }( A_0, A_1, \cdots, A_n)
 \e^{{ \l_j}}+ \sum_{k=0}^\infty \beta_k ( A_0, A_1, \cdots, A_n)
\e^k\log\ \e\cr 
&\qquad +
 \sum_{k=0}^\infty \gamma_k ( A_0, A_1, \cdots, A_n) \e^k,}$$
{\it where $\l_j=(j-   a   - \dim\ M)/ q$ with $  a \equiv 
 \ord\ A_0+ \cdots + \ord\ A_n$, $q\equiv  \ord\ Q$, and }\hfill\break
$\alpha_j ( A_0, A_1, \cdots, A_n),
\beta_k ( A_0, A_1, \cdots, A_n),
\gamma_k ( A_0, A_1, \cdots, A_n)
\in
\C$.
\item{(ii)} {\it For $j \in \N$ with  ${\rm Re}(\l_j)=
( j - a - \dim\ M)/ q <0$,
  there is a
 multi-index $(N_1, \cdots, N_n)\in \N^n$ such that}
 $$\eqalign{&\alpha_j( A_0, A_1, \cdots, A_n) = \cr
&\qquad  \sum_{j_1=  0}^{N_1} \cdots \sum_{j_n= 0}^{N_n}(-1)^{j_1+ \cdots j_n}\left[\int_{
0}^1 \cdots
\int_{ 0}^1
\left(
 \sigma_1^{j_1}
 \cdots  \sigma_{n }^{j_n}{\Gamma({- j + a +\dim\ M\over q}+
(j_1 + \cdots + j_n)
\over q \cdot j_1! j_2! \cdots j_n!} \right) d\sigma_1 \cdots
d \sigma_n\right]\cdot \cr
& \qquad\cdot
\res\left( A_0  [ A_1]_Q^{j_1} \cdots [ A_n]_Q^{j_n} Q^{{j-a - \dim\ M\over q}
-(j_1+
\cdots+ j_n) } \right),\cr}$$  
{\it where $\res$ denotes the Wodzicki
residue.}   

\bigskip

\noindent
{\bf Remark:}
%Remark (i) seems to be shown in Lemma C.3 below.
% (i) This generalizes the following known result:  
%for $n=0$, if  $a_0\equiv  \ ord\ A_0$ satisfies
%$( j - a_0 - \dim\ M/ q)\in \C-\{0, -1, -2,\cdots \}$,
%then
%$$\alpha_{j }( A_0 )={   \Gamma({ -j + a + \dim\ M\over q}) \over q}
%res( A_0  Q^{j-a - \dim\ M\over q}
%  ).$$ 
Assuming ${\rm Re}(\l_j)<0$ ensures that
  $\Gamma\left({- j + a +dim M\over q}+
(j_1 + \cdots + j_n)\right) $ is well defined for all \hfill\break
$(j_1, \cdots, j_n ) \in \N^n$.
  \medskip 

\noindent
 The proof of the theorem depends on a 
 lemma  whose proof we include for completeness.

\bigskip

\noindent\hang {\bf Lemma B.5:} {\it Let
   $A\in PDO(M,E)$ and  $Q\in Ell^+_{ord >0} (M, E)$. There is the
asymptotic expansion as $\e\to 0$}
$$\tr(A e^{-\e Q})\sim   \sum_{j=0}^\infty
 \alpha_{j }(A)
 \e^{{ \l_j}}+ \sum_{k=0}^\infty \beta_k (A)   \e^k\log\ \e+\sum_{k=0}^\infty
\gamma_k (A)\e^k,\eqno({\rm B}.3)$$
{\it with $\alpha_j(A),\ \beta_k(A),\ \gamma_k (A)\in \C$,
$a= \ord(A),\ q = \ord (Q)$, and
  $\l_j = (j-a-\dim\ M) / q$.
For  $j$ with ${\rm Re}(\l_j) <0$ (e.g.~$\ord\ A\not\in \Z$), and for
$k\in \N$, we have}
   $$\eqalign{   \alpha_{j}(A) &= {  \Gamma( {-j+ a + \dim\ M \over q})\over q}
\res(AQ^{j- a - \dim\ M
\over q}),\cr
\beta_k (A) &= (-1)^k{ q\cdot\res(A Q^{-k})  \over (k-1)!}.\cr}$$ 

\bigskip

\noindent {\sc Proof:}
For $s \in \C-\{0, -1, -2, \cdots\}$, we have
$$\eqalign{ q^{-1} \res (AQ^{-s})&= \res_{z=0}   \tr( A Q^{-(z+s)})\cr
&= \res_{z=0} \left( {1 \over \Gamma (s+z)} \int_0^\infty
t^{s +z-1}\tr(A e^{-t Q})dt\right)\cr
 &= \res_{z=0}\left(   {1\over \Gamma (s)} \int_0^1 t^{s+z-1} \tr(A e^{-t Q})dt
 + {1\over \Gamma (s)} \int_1^\infty t^{s +z-1}\tr(A e^{-t Q})dt \right)\cr
&=  \sum_j {\alpha_j(A) \over \Gamma(s)}
\res_{z=0} \left(  \int_0^1 t^{z+ \l_j  +s-1}  dt \right)\cr
&\qquad +  \sum_{k=0}^\infty {\beta_k(A) \over \Gamma(s)} \res_{z=0}
\left( \int_0^1 t^{k+   z+ s-1} \log t\   dt \right)\cr
&= \Gamma(s)^{-1}  \sum_j \alpha_j (A)\res_{z=0}
    \left[{t^{z+  \l_j +s } \over z+\l_j +s}\right]_0^1 \cr
&= \Gamma(s)^{-1} \alpha_{- q \cdot s+ a + \dim\ M}(A)
    \cr}\eqno({\rm B}.4)$$
since $\l_j=-s$ iff $j=-qs + a + \dim\ M$. Notice that $s\in -\N$ iff
$\l_j= (j-a - \dim\ M) / q \in \N$,  which does not occur if ${\rm Re}
(\l_j) <0$.
In this computation 
we use the fact that the terms in (B.3) containing logarithmic
    divergences in $\e$ or having integral powers of $\e$ do not
    contribute to the residue at
$z=0$. Similarly, the $\int_1^\infty$ term in (B.4) does not
    contribute to the residue.

\smallskip 
For $s=-l$, $l \in \N$, using 
$\Gamma(z)=  z (z-1) \cdots (z-k+1) \Gamma(z-k)$ and 
$\Gamma(z)^{-1} \sim z$ as $z\to 0$, we  find
$$\eqalign{ q^{-1} \res(A Q^{-l})&= \res_{z=0} \left( {1\over \Gamma(z-l)}
\int_0^\infty t^{z-(l+1)} \tr(A e^{-t Q}) dt \right) \cr
&= \res_{z=0} \left( {z (z-1) \cdots (z-l +1)\over \Gamma(z )} \int_0^\infty
t^{z-(l+1)} \tr(A e^{-t Q}) dt \right) \cr
&= -\sum_{k=0}^\infty \beta_k \res_{z=0} \left( z^2 (z-1) \cdots (z-l+1)
{t^{ z-l +k} \over (z-l+k)^2}\right)\cr
&= (-1)^l (l-1)! \beta_l .\cr}$$
${}$\hfill$\bullet$

\medskip \noindent { \bf Proof of the Theorem:}
(i) The operator
$A_0
  [ A_1]^{j_1}_Q  \cdots  [ A_n]^{j_n}_Q $ is  a PDO of
order at most
$a+ (j_1+ \cdots + j_n )q$, so 
$\tr(A_0
[ A_1]^{j_1}_Q  \cdots  [ A_n]^{j_n}_Q e^{- \e Q})$ has an asymptotic expansion
as in  (B.3)  with $\l_j=[(j- a- \dim\ M)/ q]
-(  j_1+ \cdots + j_n)  $. By Prop.~C.1, 
$\langle A_0, A_1, \cdots, A_n\rangle_{  \e, n, Q}$ has an asymptotic
expansion as in
(B.3) with $\l_j= (j- a- \dim\ M)/ q $. Let  $\tilde
\alpha_j (A_0,
  A_1, \cdots, A_n)$ be the coefficient of 
 $\e^{\l_j}$  in the asymptotic expansion of $\tr(A_0
  [ A_1]^{j_1}_Q  \cdots  [ A_n]^{j_n}_Q e^{- \e Q})$ 
with $\l_j=[(j- a- \dim\ M)/ q]
-(  j_1+ \cdots + j_n). $

(ii) By Lemma  B.3, if
${\rm Re}( \l_j)- (j_1+ \cdots + j_n) <0$ 
(e.g.~if ${\rm Re}( \l_j)<0$),
 then 
$$\eqalign{ \tilde \alpha_j (A_0,
    A_1, \cdots  , A_n,) &=
  {  \Gamma(-\l_j+  (j_1+ \cdots + j_n))\over q} \cr
&\qquad \cdot
\res(A_0
  [ A_1]^{j_1}_Q  \cdots
 [ A_n]^{j_n}_Q Q^{{j- a - \dim\ M \over q}- (j_1+ \cdots + j_n)}  ).\cr}
$$
Part (ii) of the theorem  follows.\hfill $\bullet$
    \vfill \eject  

\centerline{\bf References}
\sp \noi

 \item{[AJPS]} S. Albeverio, J. Jost, S. Paycha, S. Scarlatti,
{\it  A Mathematical Introduction to String Theory},
LMS Lecture Note Series {\bf 225}, Cambridge, UK, Cambridge University
Press, 1997.\medskip

\item{[AP]} M. Arnaudon, S. Paycha, ``Regularisable and minimal orbits
for group actions in infinite dimensions,''     {\it Commun. Math. Phys.}
{\bf 191} (1998),
641-662.\medskip

 \item{[B]} J.-M. Bismut, ``Localization formulae, super connections and
the index theorem for families,'' {\it Commun. Math. Phys.} {\bf 103} (1986),
127-166.\medskip

\item{[BF]} J.-M. Bismut, D. Freed, ``The analysis of elliptic 
families I,''
{\it Commun. Math. Phys.} {\bf 106} (1986), 159-176.\medskip

 \item{[BGV]} N. Berline, E. Getzler, M. Vergne, { \it Heat Kernels and
Dirac Operators}, Berlin, Springer-Verlag, 1991.\medskip

 \item{[CDMP]} A. Cardona, C. Ducourtioux, J. P. Magnot, S. Paycha, ``Weighted
traces on   algebras  of pseudodifferential operators and geometry of loop
groups'', preprint, 1999.\medskip

\item{[CM]} A. Connes, H. Moscovici, ``The local index formula in 
non-commutative geoemetry,'' {\it GAFA} {\bf 5}
(1995), 174-243.
\medskip

\item{[F]} D. Freed, ``The geometry of loop groups'', {\it Journal Diff. Geom.}
{\bf 28} (1988),
  223--276.\medskip

\item{[G]} P. B. Gilkey, {\it Invariance Theory, the Heat Equation,
  and the Aityah-Singer Index Theorem}, Houston, Publish or Perish, 1984.

\medskip
\item{[JLO]} A. Jaffe, A. Lesniewski, K. Osterwalder, ``Quantum K-theory,
I. The Chern character,'' {\it Commun Math. Phys.} {\bf 118} (1988),
1-14.\medskip 

\item{[K]} C. Kassel, ``Le r\'esidu non commutatif [d'apr\`es
Wodzicki],'' S\'eminaire Bourbaki {\bf 708}, {\it Asterisque} (1989),
199-229. \medskip

\item{[KV1]} M. Kontsevich, S. Vishik, ``Determinants of elliptic
pseudodifferential
operators,'' Geometry and Functional Analysis, Max Planck Preprint,
1994.\medskip

\item{[KV2]} M. Kontsevich, S. Vishik, ``Geometry of determinants of
elliptic operators'' in {\it Functional Analysis on the Eve of the
21st Century, Vol. I, }
 (ed. S. Gindikin, J. Lepowski, R. L. Wilson)  Progress in
Mathematics, Birkh\"auser, 1994, 173-197.

\item{[L]} M. Lesch, ``On the non-commutative
residue for pseudo-differential operators with
log-polyhomogeneous symbols,'' {\it Annals of
Global Anal. and Geom.} {\bf 17} (1998), 151-187.\medskip
  
\item{ [MRT]} Y. Maeda, S. Rosenberg, P. Tondeur,
``The mean curvature of gauge orbits,'' in {\it Global Analysis and Modern
Mathematics} (ed.
 K. Uhlenbeck), Houston, Publish or Perish, 1994, 171-220.\medskip

\item{} Y. Maeda, S. Rosenberg, P. Tondeur,  ``Minimal orbits of metrics,''
{\it J.
Geom. and Phys.} {\bf 23} (1997),  314-349.\medskip

\item{[Pa]} S. Paycha, ``Regularized  pseudo-traces as a looking glass
into infinite dimensional
geometry: a proposal
to extend some concepts of Riemannian geometry to infinite dimensions,''
preprint, 2000.\medskip

 \item{[PR]} S. Paycha, S. Rosenberg, 
``About infinite dimensional group actions and determinant bundles,'' in
{\it Analysis on Infinite-Dimensional Lie
Groups and Algebras}, (eds. H. Heyer, J. Marion)
Singapore, World Scientific, 1998, 355-367.\medskip

\item{[Q1]} D. Quillen, ``Determinants of  Cauchy-Riemmann operators
over a Riemann surface,'' {\it Funktsional Anal. i Prilozhen.} {\bf
19} (1985), 37-41.

\medskip
\item{[Q2]} D. Quillen, ``Superconnections and the Chern character,''
{\it Topology} {\bf 24} (1985), 89-95

\medskip
 \item{[T]}    T. Tromba, {\it Teichmü\"uller
Theory in Riemannian Geometry}, Boston, Birkh\"auser Verlag, 1992.\medskip

\item{[W]} M. Wodzicki, ``Non-commutative residue'',  Lecture Notes in
Mathematics {\bf 1289}, Berlin, Springer-Verlag, 1987.\medskip
 \bye